\documentstyle[12pt,a4wide,leqno,amssymb,bezier]{article}

\begin{document}

\title{\bf Polynomial growth of the derivative for diffeomorphisms on
tori}
\author{Krzysztof Fr\c{a}czek}
\maketitle
\renewcommand{\thefootnote}{}
\footnote{2000
{\em Mathematics Subject Classification}: 37A05, 37C05, 37C40.}
\footnote{Research partly supported by KBN grant 5 P03A 027 21(2001).}
\newcommand{\tor}{{\Bbb{T}}}
\newcommand{\ep}{\varepsilon}
\newcommand{\pf}{{\bf Proof. }}
\newtheorem{prop}{\indent Proposition}
\newtheorem{theo}{\indent Theorem}
\newtheorem{lem}[theo]{\indent Lemma}
\newtheorem{cor}{\indent Corollary}
\newtheorem{df}{\indent Definition}

\begin{abstract}
We consider area--preserving diffeomorphisms on tori
with zero entropy. We classify ergodic area--preserving
diffeomorphisms of the $3$--torus for which the sequence
$\{Df^n\}_{n\in{\Bbb N}}$ has polynomial growth. Roughly speaking, the
main theorem says that every ergodic area--preserving
$C^2$--diffeomorphism with polynomial uniform growth of the derivative is
$C^2$--conjugate to a $2$--steps skew product of the form
\[\tor^3\ni(x_1,x_2,x_3)\mapsto
(x_1+\alpha,\ep x_2+\beta(x_1),x_3+\gamma(x_1,x_2))\in\tor^3,\]
where $\ep=\pm 1$.
We also indicate why there is no $4$--dimensional analogue of
the above result. Random diffeomorphisms on the $2$--torus are studied as well.
\end{abstract}

\section{Introduction}\indent

Let $M$ be a compact Riemannian smooth manifold and let $\mu$ be a probability
Borel measure with supp$(\mu)=M$. Let  $f:(M,\mu)\rightarrow(M,\mu)$ be a
smooth measure--preserving diffeomorphism.
An important question of smooth ergodic theory is: what is
the relation between
asymptotic properties of the sequence $\{Df^n\}_{n\in{\Bbb N}}$ and dynamical
properties of the dynamical system  $f:(M,\mu)\rightarrow(M,\mu)$.
There are results well describing this relation in the case where $M$ is the
torus. For example, if $f$ is homotopic to
the identity, the rotation vector of $f$ is ergodic and the sequence
$\{Df^n\}_{n\in{\Bbb N}}$ is uniformly bounded, then
$f$ is $C^0$--conjugate to an ergodic rotation (see \cite{Her} p.181).
Moreover, if $\{Df^n\}_{n\in{\Bbb N}}$ is bounded in
the $C^{r}$--norm ($r\in{\Bbb N}\cup\{\infty\}$), then $f$ and the ergodic
rotation are $C^{r}$--conjugated (see \cite{Her} p.182).
On the other hand, if $\{Df^n\}_{n\in{\Bbb N}}$ has "exponential growth",
precisely if $f$ is an Anosov diffeomorphism, then $f$ is $C^0$--conjugate
to an algebraic automorphism of the torus (see \cite{Man}).

A natural question is: what can happen between the
above extreme cases?
The aim of this paper is to classify measure--preserving
diffeomorphisms $f$ of tori for which the sequence
$\{Df^n\}_{n\in{\Bbb N}}$ has polynomial growth.
 One definition of the
polynomial growth of the derivative is presented in \cite{Fr}.
In the above--mentioned paper, the following result is proved.

\begin{prop}\label{wstep}
Let $f:\tor^2\rightarrow\tor^2$ be an ergodic area--preserving
$C^2$--diffeomorphism.
If the sequence $\{n^{-\tau}Df^n\}_{n\in{\Bbb N}}$ converges
a.e.\ ($\tau>0$) to a measurable nonzero function, then $\tau=1$ and $f$ is
algebraically conjugate (i.e.\ via a group automorphism) to the skew
product of an irrational rotation on the circle and a circle
cocycle with nonzero topological degree.
\end{prop}

Moreover, in~\cite{Fr1}, the
author has shown that  if $f:\tor^2\rightarrow\tor^2$ is an ergodic
area--preserving $C^3$--diffeomorphism for which the sequence
$\{n^{-1}Df^n\}_{n\in{\Bbb N}}$ is uniformly
separated from $0$ and $\infty$ and it is bounded in the $C^2$--norm, then
$f$ is also algebraically conjugate to the skew
product of an irrational rotation  on the circle and a circle
cocycle with nonzero topological degree.

Next very interesting result indicating how the type of growth of the
derivative of a diffeomorphism reflects its important dynamical features
(see~\cite{Po-Si} and \cite{Po}) says that every symplectic diffeomorphism
$f:\tor^2\rightarrow\tor^2$
homotopic to the identity with  a fixed point either is equal to the
identity map or there exists $c>0$ such that
\[\max(\|Df^n\|_{sup},\|Df^{-n}\|_{sup})\geq cn\]
for any natural $n$.

In this paper some analogous of Proposition~\ref{wstep} are studied.
In Section~\ref{random} we discuss some random versions of
Proposition~\ref{wstep}. In Section~\ref{3d} we
classify area--preserving
ergodic $C^2$--diffeomorphisms of polynomial uniform growth of the derivative
on the $3$--torus, i.e.\ diffeomorphisms for which the sequence
$\{n^{-\tau}Df^n\}_{n\in{\Bbb{N}}}$ converges uniformly to a non--zero
function. It is shown that
if the limit function is of class $C^1$, then $\tau$ equals either $1$ or
$2$ and the diffeomorphism is $C^2$--conjugate to a $2$--steps skew product.
In Section~\ref{4d} we indicate why there is no $4$--dimensional analogue of
the above result.

\section{Random diffeomorphism on the $2$--torus}\indent\label{random}

Throughout this section we will consider smooth random dynamical
systems over an abstract dynamical system
$(\Omega,{\cal{F}},P,T)$, where $(\Omega,{\cal{F}},P)$ is a
Lebesgue space and
$T:(\Omega,{\cal{F}},P)\rightarrow(\Omega,{\cal{F}},P)$ is an
ergodic measure--preserving automorphism. As a phase space for
smooth random diffeomorphisms we will consider a compact
Riemannian $C^{\infty}$--manifold $M$ equipped with its Borel
$\sigma$--algebra ${\cal{B}}$. A measurable map $f$
\[{\Bbb{Z}}\times\Omega\times M\ni(n,\omega,x)\longmapsto f^n_{\omega}x\in M\]
satisfying for $P$--a.e.\ $\omega\in\Omega$ the following conditions
\begin{itemize}
\item $f^0_{\omega}=\mbox{ Id}_M$,
$f^{m+n}_{\omega}=f^{m}_{T^n\omega}\circ f^{n}_{\omega}$ for all
$m,n\in{\Bbb{Z}}$,
\item $f^n_{\omega}:M\rightarrow M$ is a smooth function for all
$n\in{\Bbb{Z}}$,
\end{itemize}
is called a {\em smooth random dynamical system} (RDS). Of course,
the smooth RDS is generated by the random diffeomorphism
$f_{\omega}=f^1_{\omega}$ in the sense that
\[f^n_{\omega}=\left\{
\begin{array}{rcl}
f_{T^{n-1}\omega}\circ\ldots\circ f_{T\omega}\circ f_{\omega} & \mbox{for}
& n>0\\
\mbox{ Id}_M & \mbox{for} & n=0\\
f_{T^{n}\omega}^{-1}\circ f_{T^{n+1}\omega}^{-1}\circ\ldots\circ
f_{T^{-1}\omega}^{-1} & \mbox{for} & n<0.
\end{array}
\right.\]
Consider the skew--product transformation
$T_f:(\Omega\times M,{\cal{F}}\otimes{\cal{B}})\rightarrow
(\Omega\times M,{\cal{F}}\otimes{\cal{B}})$
induced naturally by $f$
\[T_f(\omega,x)=(T\omega,f_{\omega}x).\]
Then $T^{n}_f(\omega,x)=(T^{n}\omega,f^{n}_{\omega}x)$ for all
$n\in{\Bbb{Z}}$. We call a probability measure $\mu$ on
$(\Omega\times M,{\cal{F}}\otimes{\cal{B}})$ {\em $f$--invariant},
if $\mu$ is invariant  under $T_f$ and has marginal $P$ on
$\Omega$. Such measures can also be characterized in terms of
their disintegrations $\mu_{\omega}$, $\omega\in\Omega$ by
$f_{\omega}\mu_{\omega}=\mu_{T\omega}$ $P$--a.e. A $f$--invariant
measure $\mu$ is said to be {\em ergodic} if
$T_f:(\Omega\times
M,{\cal{F}}\otimes{\cal{B}},\mu)\rightarrow(\Omega\times
M,{\cal{F}}\otimes{\cal{B}},\mu)$ is ergodic. We say that a $f$--invariant
measure $\mu$ has full support, if supp$(\mu_{\omega})=M$
for $P$--a.e.\ $\omega\in\Omega$.

In this section we will deal with almost everywhere diffentiable and $C^r$
measure--preserving random dynamical systems with polynomial growth of the
derivative.
Suppose that $f:{\Bbb{Z}}\times\Omega\times M\rightarrow M$ is a $C^0$  RDS
and $\mu$ is a $f$--invariant measure on $\Omega\times M$.
The RDS $f$ is called $\mu$--almost everywhere diffentiable if for every
integer $n$ and for $\mu$--a.e.\ $(\omega,x)\in\Omega\times M$ there
exists the derivative $Df^{n}_{\omega}(x):T_xM\rightarrow
T_{f^n_{\omega}}M$ and
\[\int_M\|Df^{n}_{\omega}(x)\|_{n,\omega,x}d\mu_{\omega}(x)<\infty\]
for every $n\in{\Bbb{Z}}$ and $P$--a.e.\ $\omega\in\Omega$, where
$\|\,\cdot\,\|_{n,\omega,x}$ is the operator norm in
${\cal{L}}(T_xM,T_{f^n_{\omega}x}M)$.

In the paper we will look more closely at RDS on tori.
Let $d$ be a natural number. By $\tor^d$ we will mean the $d$--dimensional
torus
$\{(z_1,\ldots,z_d)\in{\Bbb{C}^d}:|z_1|=\ldots=|z_d|=1\}$ which
most often will be treated as the quotient group
${\Bbb{R}}^d/{\Bbb{Z}}^d$; $\lambda^{\otimes d}$ will denote Lebesgue
measure on $\tor^d$.
We will identify functions on $\tor^d$ with ${\Bbb{Z}}^d$--periodic
functions (i.e.\ periodic of period 1 in each coordinates) on
${\Bbb{R}}^d$.
Let $f:\tor^d\rightarrow\tor^d$ be a
smooth diffeomorphism. We will identify $f$ with a diffeomorphism
$f:{\Bbb{R}}^d\rightarrow{\Bbb{R}}^d$ such that
\[f(x_1,\ldots,x_j+1,\ldots,x_d)=f(x_1,\ldots,x_d)+(a_{1j},\ldots,a_{dj})\]
for every $(x_1,\ldots,x_d)\in{\Bbb{R}}^d$, where
$A=[a_{ij}]_{1\leq i,j\leq d}\in GL_d({\Bbb{Z}})$. We call $A$ the {\em linear
part} of the diffeomorphism $f$.
Then there exist smooth functions
$\tilde{f}_i:\tor^d\rightarrow{\Bbb{R}}$  such that
\[f_i(x_1,\ldots,x_d)=\sum_{j=1}^da_{ij}x_j+\tilde{f}_i(x_1,\ldots,x_d),\]
where $f_i:{\Bbb{R}}^d\rightarrow{\Bbb{R}}$ is the $i$--th
coordinate functions of $f$.

\begin{df}
{\em We say that a $\mu$--almost everywhere diffentiable RDS $f$ on $\tor^d$
over  $(\Omega,{\cal{F}},P,T)$
has} $\tau$--polyno\-mial ($\tau>0$) growth of the derivative {\em if
\[\frac{1}{n^{\tau}}Df^n_{\omega}(x)\rightarrow g(\omega,x) \mbox{ for
$\mu$--a.e. } (\omega,x)\in \Omega\times \tor^d,\]
where $g:\Omega\times \tor^d\rightarrow M_d({\Bbb{R}})$ is $\mu$ non--zero,
i.e.\ there exists a set $A\in{\cal{F}}\otimes{\cal{B}}$ such that
$\mu(A)>0$ and $g(x)\neq 0$ for all $x\in A$. Moreover, if additionally
$Df^n${ belongs to $L^1((\Omega\times\tor^d,\mu),M_d({\Bbb{R}}))$ for all
} $n\in{\Bbb{N}}$ and the sequence
$\{n^{-\tau}Df^n\}$ converges in $L^1((\Omega\times\tor^d,\mu),M_d({\Bbb{R}}))$
then we say that $f$ has } $\tau$--polynomial L$^1$--growth of the derivative.
\end{df}

We now give the exmple of an ergodic RDS on $\tor^2$ with linear
$L^1$--growth of the derivative.
For  given $\tau:X\rightarrow X$, $\varphi:X\rightarrow\tor$ and $n\in{\Bbb{N}}$ let
\[\varphi^{(n)}=\varphi(x)+\varphi(\tau x)+\ldots+\varphi(\tau^{n-1}x).\]
Let us consider an almost everywhere diffentiable RDS $f$ on $\tor^2$ over
$(\Omega,{\cal{F}},P,T)$ (called the random Anzai skew product) of
the form
\[f_{\omega}(x_1,x_2)=(x_1+\alpha(\omega),x_2+\varphi(\omega,x_1)),\]
where the skew product $T_\alpha:
(\Omega\times\tor,P\otimes\lambda)\rightarrow
(\Omega\times\tor,P\otimes\lambda)$,
$T_{\alpha}(\omega,x)=(T\omega,x+\alpha(\omega))$ is ergodic and
$\varphi:\Omega\times\tor\rightarrow\tor$ is an absolutely
continuous random mapping of the circle such that $D\varphi\in
L^1(\Omega\times\tor,P\otimes\lambda)$ and
$\int_{\Omega}d(\varphi_{\omega})dP\omega\neq 0$.
Then the product measure $P\otimes\lambda^{\otimes 2}$ is $f$--invariant.
The following lemma is a little bit more general version of Lemma 3 in
\cite{Iw-Le-Ru}.

\begin{lem}\label{skew}
The RDS $f$ is ergodic and has linear L$^1$--growth of the derivative.
\end{lem}

\pf First note that
\[f^n_{\omega}(x_1,x_2)=(x_1+\alpha^{(n)}(\omega),x_2+\varphi^{(n)}(\omega,x_1))\]
for all $n\in{\Bbb{N}}$. Therefore
\[\frac{1}{n}Df^n_{\omega}(x_1,x_2)=\left[
\begin{array}{cc}
1/n & 0 \\
1/n\sum_{k=0}^{n-1}D\varphi(T_{\alpha}^k(\omega,x_1)) & 1/n
\end{array}
\right].\]
By the ergodicity of $T_{\alpha}$,
\[\frac{1}{n}\sum_{k=0}^{n-1}D\varphi(T^k_{\alpha}(\omega,x))\rightarrow
\int_{\Omega}\int_{\tor}D\varphi_{\omega}(y)dydP\omega=
\int_{\Omega}d(\varphi_{\omega})dP\omega\neq 0\]
for $P\otimes\lambda$--a.e.\ $(\omega,x)\in\Omega\times\tor$ and
in the L$^1$--norm, which implies linear L$^1$--growth of the derivatives
of $f$.

To proof the ergodicity of $f$, we consider the family of unitary
operators $U_m:L^2(\Omega\times\tor,P\otimes\lambda)\rightarrow
L^2(\Omega\times\tor,P\otimes\lambda)$, $m\in{\Bbb{Z}}$ given by
$U_mg(\omega,x)=e^{2\pi
im\varphi(\omega,x)}g(T\omega,x+\alpha(\omega))$. We now show that
\begin{equation}\label{uny}
\langle U^n_mg,g\rangle=\int_{\Omega\times\tor}
e^{2\pi im\varphi^{(n)}(\omega,x)}gT_{\alpha}^n(\omega,x)
\bar{g}(\omega,x)dP\omega dx\rightarrow 0\mbox{ as }n\rightarrow\infty
\end{equation}
for all $g\in L^2(\Omega\times\tor,P\otimes\lambda)$ and
$m\in{\Bbb{Z}}\setminus\{0\}$.
Let $\Lambda$ denote the set of all $g\in
L^2(\Omega\times\tor,P\otimes\lambda)$ satisfying (\ref{uny}). It is easy
to check that $\Lambda$ is a closed linear subspace of
$L^2(\Omega\times\tor,P\otimes\lambda)$. Therefore it suffices to show
(\ref{uny}) for all functions of the form $g(\omega,x)=h(\omega)e^{2\pi ikx}$,
where $h\in L^{\infty}(\Omega,P)$ and $k\in{\Bbb{Z}}$. Then
\begin{eqnarray*}
|\langle U^n_mg,g\rangle| & = & |\int_{\Omega}h(T^n\omega)
\bar{h}(\omega)e^{2\pi ik\alpha^{(n)}(\omega)}(\int_{\tor}
e^{2\pi im\varphi^{(n)}(\omega,x)}dx)dP\omega |\\
 & \leq & \|h\|_{L^{\infty}}^2\int_{\Omega}|\int_{\tor}
e^{2\pi im\varphi^{(n)}(\omega,x)}dx|dP\omega.
\end{eqnarray*}
Let $\tilde{\varphi}:\Omega\times\tor\rightarrow{\Bbb{R}}$ be an
absolutely continuous random function such that
$\varphi(\omega,x)=\tilde{\varphi}(\omega,x)+d(\varphi_{\omega})x$.
Without loss of generality we can assume that
$\int_{\Omega}d(\varphi_{\omega})dP\omega=a>0$. For any natural $n$ let
\[A_n=\{\omega\in\Omega:d(\varphi_{\omega})^{(n)}/n>a/2\}.\]
By the ergodicity of $T$, $P(\Omega\setminus A_n)\rightarrow 0$ as
$n\rightarrow\infty$. Applying integration by parts we obtain
\begin{eqnarray*}
\frac{1}{\|h\|_{L^{\infty}}^2}|\langle U^n_mg,g\rangle| & \leq & P(\Omega\setminus
A_n)+\int_{A_n}|\int_{\tor}e^{2\pi im\tilde{\varphi}^{(n)}(\omega,x)}
d\frac{e^{2\pi im{d(\varphi_{\omega})}^{(n)}x}}
{2\pi im{d(\varphi_{\omega})}^{(n)}}|dP\omega \\
 & \leq & P(\Omega\setminus A_n)+\frac{1}{\pi |m|an}
\int_{A_n}|\int_{\tor}e^{2\pi im{d(\varphi_{\omega})}^{(n)}x}
de^{2\pi im\tilde{\varphi}^{(n)}(\omega,x)}|dP\omega \\
 & \leq & P(\Omega\setminus A_n)+\frac{2}{\pi an}
\int_{A_n}|\int_{\tor}D\tilde{\varphi}^{(n)}(\omega,x)dx|dP\omega \\
 & \leq & P(\Omega\setminus A_n)+\frac{2}{\pi a}
\int_{\Omega\times\tor}|D\tilde{\varphi}^{(n)}(\omega,x)/n|dP\omega dx.
\end{eqnarray*}
As $\int_{\Omega\times\tor}D\tilde{\varphi}(\omega,x)dP\omega dx=0$, applying
Birkhoff's Ergodic Theorem for $T_{\alpha}$ we conclude that
$\int_{\Omega\times\tor}|D\tilde{\varphi}^{(n)}(\omega,x)/n|dP\omega dx$
tends to zero, which proves our claim.

Now suppose, contrary to our assertion, that $f$ is not ergodic.
Since the skew product $T_{\alpha}$ is ergodic, there exists a
measurable function $g:\Omega\times\tor\rightarrow\tor$ and
$m\in{\Bbb{Z}}\setminus\{0\}$ such that
\[e^{2\pi
im\varphi(\omega,x)}=g(\omega,x)\bar{g}(T_{\alpha}(\omega,x)).\]
Then $\langle U^n_mg,g\rangle=1$ for all $n\in{\Bbb{N}}$, contrary
to (\ref{uny}). $\Box$

\vspace{2ex}

The aim of this section is to classify $C^r$
random dynamical systems on the $2$--torus that have polynomial
($L^1$) growth of the derivative and are ergodic with respect to an
invariant measure of full support. We say that two random dynamical
systems
$f$ and $g$ on $\tor^2$ over $(\Omega,{\cal{F}},P,T)$ are
algebraically conjugated if there exists a group automorphism
$A:\tor^2\rightarrow\tor^2$ such that
$f_{\omega}\circ A=A\circ g_{\omega}$ for $P$--a.e.\ $\omega
\in\Omega$. We prove the following theorems.

\begin{theo}\label{rtg}
Let $f$ be a $C^r$ random dynamical system on $\tor^2$ over
$(\Omega,{\cal{F}},P,T)$. Let $\mu$ be an $f$--invariant ergodic
 measure of full support on $\Omega\times\tor^2$. Suppose that $f$ has
$\tau$--polynomial growth of the derivative. Then $\tau\geq 1$
and  $f$ is algebraically conjugate to a random skew product of
the form
\[\hat{f}_{\omega}(x_1,x_2)=(F_{\omega}(x_1),x_2+\varphi_{\omega}(x_1)),\]
where $F:\Omega\times\tor\rightarrow\tor$ is a $C^r$ random
diffeomorphism of the circle. Moreover, if $f$ preserves
orientation, then there exist a random homeomorphism of the
circle $\xi:\Omega\times\tor\rightarrow\tor$ and  a measurable
function $\alpha:\Omega\rightarrow\tor$ such that
\[\xi_{T\omega}\circ F_{\omega}(x) =
\xi_{\omega}(x)+\alpha_{\omega}\mbox{ $P$--a.e.}\] and
consequently $f$ is topologically conjugate to the random skew
product
\[\tor^2\ni(x_1,x_2)\longmapsto(x_1+\alpha_{\omega},
x_2+\varphi_{\omega}\circ\xi^{-1}_{\omega}(x_1))\in\tor^2.\]
\end{theo}
\begin{theo}\label{rtg1}
Under the hypothesis of Theorem~\ref{rtg}, if moreover $f$ has
$\tau$--polynomial $L^1$--growth of the derivative and $\mu$ is
equivalent to the measure $P\otimes\lambda^{\otimes 2}$ with
$d\mu/d(P\otimes\lambda^{\otimes 2})$, $d(P\otimes\lambda^{\otimes 2})/d\mu\in
L^{\infty}(\Omega\times\tor^2)$ , then
\begin{itemize}
\item $\tau=1$,
\item there exist a Lipschitz random diffeomorphism of the
circle $\xi:\Omega\times\tor\rightarrow\tor$ with $D\xi,D\xi^{-1}\in
L^{\infty}(\Omega\times\tor,P\otimes\lambda)$ and a measurable
function $\alpha:\Omega\rightarrow\tor$ such that
\[\xi_{T\omega}\circ F_{\omega}(x) =
\xi_{\omega}(x)+\alpha_{\omega}\mbox{ $P$--a.e.\ and}\]
\item $\int_{\Omega}d(\varphi_{\omega}\circ\xi^{-1}_{\omega})dP\omega\neq
0$.
\end{itemize}
\end{theo}

For convenience, the proofs of the theorems are divided into a sequence of
lemmas.
Let $f$ be a $C^1$ random dynamical system on $\tor^d$ over
$(\Omega,{\cal{F}},P,T)$. Let $\mu$ be a $f$--invariant ergodic
 measure of full support on $\Omega\times\tor^d$. Suppose that $f$ has
$\tau$--polynomial growth of the derivative. Let
$g:\Omega\times \tor^d\rightarrow M_{d}({\Bbb{R}})$ denote the limit
of the sequence $\{n^{-\tau}Df^{n}\}$.

\begin{lem}\label{lemg}
For $\mu$--a.e.\ $(\omega,x)\in\Omega\times \tor^d$ and all integer $n$ we
have
\begin{equation}\label{dwa}
g(\omega,x)\neq 0,\;\;\;\;g(\omega,x)^2=0\mbox{ and}
\end{equation}
\begin{equation}\label{trzy}
g(\omega,x)=g(T^n\omega,f^{n}_{\omega}x)Df^n_{\omega}(x).
\end{equation}
For $\mu\otimes\mu$--a.e.\ $(\omega,x,\upsilon,y)\in
\Omega\times \tor^d\times\Omega\times \tor^d$ we have
\begin{equation}\label{cztery}
g(\omega,x)g(\upsilon,y)=0\mbox{ and }g(\omega,x)=Df_{\upsilon}(y)\,g(\omega,x).
\end{equation}
\end{lem}

\pf
Let $A\subset\Omega\times\tor^d$ be a full $\mu$--measure $T_f$--invariant
subset such that  $(\omega,x)\in A$ implies
$\lim_{n\rightarrow\infty}n^{-\tau}Df_{\omega}^n(x)=g(\omega,x)$.
Assume that $(\omega,x)\in A$. Since
\[(\frac{m+n}{m})^{\tau}\frac{1}{(m+n)^{\tau}}Df_{\omega}^{m+n}({x})=
\frac{1}{m^{\tau}}Df_{T^n\omega}^m(f_{\omega}^n{x})\:
Df_{\omega}^n({x})\]
and $(T^n\omega,f_{\omega}^n{x})\in A$ for all $m,n\in{\Bbb{N}}$, letting
$m\rightarrow\infty$, we obtain
\[g(\omega,x)=g(T^n\omega,f_{\omega}^nx)\:Df_{\omega}^n({x})\;\;\mbox{
for $(\omega,x)\in A$ and $n\in\Bbb{N}$.}\]
Let $B=\{(\omega,x)\in A:g(\omega,x)\neq 0\}$. By the above, $B$ is
$T_f$--invariant. As $g$ is $\mu$
non--zero, we have $\mu(B)=1$, by the ergodicity of $T_f$.

By the Jewett--Krieger Theorem, we can assume that $\Omega$ is a compact
metric space, $T:\Omega\rightarrow\Omega$ is an uniquely ergodic
homeomorphism  and $P$ is the unique $T$--invariant measure.
Now choose a sequence $\{A_k\}_{k\in\Bbb{N}}$ of measurable subsets of $A$
such that the functions $g,Df:A_k\rightarrow
M_{d}(\Bbb{R})$ are continuous, all non-empty open subsets of
$A_k$ (in the induced topology) have positive measure and $\mu(A_k)>1-1/k$
for any natural $k$. Since the transformation
$(T_f)_{A_k}:(A_k,\mu_{A_k})\rightarrow(A_k,\mu_{A_k})$ induced by $T_f$
on $A_k$ is ergodic, for every natural $k$ we can find a measurable subset
$B_k\subset A_k$ such that for any $(\omega,x)\in B_k$ the orbit
$\{(T_f)^n_{A_k}(\omega,x)\}_{n\in\Bbb{N}}$ is dense in $A_k$ in the induced
topology and $\mu(B_k)=\mu(A_k)$.\\
Let $(\omega,x),(\upsilon,y)\in B_k$. Then there exists an increasing
sequence $\{m_i\}_{i\in\Bbb{N}}$ of natural numbers such that
$(T_f)_{A_k}^{m_i}(\omega,x)\rightarrow(\upsilon,y)$. Hence there exists
an increasing sequence $\{n_i\}_{i\in\Bbb{N}}$ of natural numbers such that
$T_f^{n_i}(\omega,x)\rightarrow(\upsilon,y)$ and $T_f^{n_i}(\omega,x)\in
A_k$ for all $i\in\Bbb{N}$. Since $g,Df:A_k\rightarrow M_{
d}(\Bbb{R})$ are continuous, we get
$g(T^{n_i}\omega,f_{\omega}^{n_i}x)\rightarrow g(\upsilon,y)$ and
$Df_{T^{n_i}\omega}(f_{\omega}^{n_i}x)\rightarrow g(\upsilon,y)$. Since
\[\frac{1}{n_i^{\tau}}g(\omega,x)=g(T^{n_i}\omega,f_{\omega}^{n_i}x)\:
\frac{1}{n_i^{\tau}}Df_{\omega}^{n_i}(x),\]
letting $i\rightarrow\infty$, we obtain $g(\upsilon,y)\:g(\omega,x)=0$.
Since
\[\frac{1}{n_i^{\tau}}Df_{\omega}^{n_i+1}({x})=
Df_{T^{n_i}\omega}(f_{\omega}^{n_i}{x})\:
\frac{1}{n_i^{\tau}}Df_{\omega}^{n_i}({x}),\]
letting $i\rightarrow\infty$, we obtain
$g(\omega,x)=Df_{\upsilon}y\:g(\omega,x)$.
Therefore
\[\mu\otimes\mu\{(\omega,x,\upsilon,y)\in
\Omega\times \tor^d\times\Omega\times \tor^d:\;
g(\upsilon,y)\:g(\omega,x)=0\})>(1-\frac{1}{k})^2,\]
\[\mu\{(\omega,x)\in
\Omega\times \tor^d:\;
g(\omega,x)^2=0\})>1-\frac{1}{k}\]
and
\[\mu\otimes\mu\{(\omega,x,\upsilon,y)\in
\Omega\times \tor^d\times\Omega\times \tor^d:\;
g(\omega,x)=Df_{\upsilon}(y)\:g(\omega,x)\})>(1-\frac{1}{k})^2\]
for any natural $k$, which proves the lemma. $\Box$

\vspace{2ex}

Let us return to the case where $d=2$.
Suppose that $A,B$ are non--zero real $2\times 2$--matrixes such that
$A^2=B^2=AB=0$.
Then (see Lemma 4 in \cite{Fr}) there exist real numbers $a,b\neq 0$ and
$c$ such that
\[\begin{array}{lcr}
A=a\left[
\begin{array}{c}
c \\
1
\end{array}\right]
\left[
\begin{array}{cc}
1 & -c
\end{array}\right]

& \mbox{ and } &
B=b\left[
\begin{array}{c}
c \\
1
\end{array}\right]
\left[
\begin{array}{cc}
1 & -c
\end{array}\right]
\end{array} \]
or
\[
\begin{array}{lcr}
A=a\left[
\begin{array}{c}
1 \\
0
\end{array}\right]
\left[
\begin{array}{cc}
0 & 1
\end{array}\right]
& \mbox{ and } &
B=b\left[
\begin{array}{c}
1 \\
0
\end{array}\right]
\left[
\begin{array}{cc}
0 & 1
\end{array}\right]
\end{array}\]
It follows that g can be represented as follows
\[
g=h \left[
\begin{array}{c}
c \\
1
\end{array}\right]
\left[
\begin{array}{cc}
1 & -c
\end{array}\right],
\]
where $h:\Omega\times\tor^2\rightarrow\Bbb{R}$ is a measurable function,
which is non--zero at $\mu$--a.e.\ point
and $c\in\Bbb{R}$.
We can omit the second case where
\[
g=h
\left[
\begin{array}{c}
1 \\
0
\end{array}\right]
\left[
\begin{array}{cc}
0 & 1
\end{array}\right],\]
because it reduces to case $c=0$ after interchanging the coordinates,
which is an algebraic isomorphism.
Then from (\ref{cztery}) we obtain
\begin{equation}\label{piec}
\left[
\begin{array}{c}
c \\
1
\end{array}\right]
=Df_{\omega}x
\left[
\begin{array}{c}
c \\
1
\end{array}\right]
\end{equation}
for $P$--a.e.\ $\omega\in\Omega$ and for all $x\in\tor^2$, because $\mu$
has full support.
>From (\ref{trzy}) we obtain
\begin{equation}\label{szesc}
h(\omega,x)\left[
\begin{array}{cc}
1 & -c
\end{array}\right]=
h(T\omega,f_{\omega}x)
\left[
\begin{array}{cc}
1 & -c
\end{array}\right]
Df_{\omega}x
\end{equation}
for $\mu$--a.e.\ $(\omega,x)\in\Omega\times\tor^2$.

\begin{lem}\label{irrat}
If $c$ is irrational, then
$f_{\omega}(x_1,x_2)=(x_1+\alpha(\omega),x_2+\gamma(\omega))$,  where
$\alpha,\gamma:\Omega\rightarrow\tor$ are measurable functions.
Consequently, the sequence $n^{-\tau}Df^n$ tends uniformly to zero.
\end{lem}

\pf
>From (\ref{piec}) we have
\[c=c\frac{\partial (f_{\omega})_1}{\partial x_1}+
\frac{\partial (f_{\omega})_1}{\partial x_2}\mbox{ and }
1=c\frac{\partial (f_{\omega})_2}{\partial x_1}+
\frac{\partial (f_{\omega})_2}{\partial x_2}\]
for $P$--a.e.\ $\omega\in\Omega$. It follows that for $i=1,2$ there exists
a $C^{r+1}$ random function
$u_i:\Omega\times{\Bbb{R}}\rightarrow{\Bbb{R}}$ such that
\[f_i(\omega,x_1,x_2)=x_i+u_i(\omega,x_1-cx_2).\]
Represent $f$ as follows
\begin{eqnarray*}
f_1(\omega,x_1,x_2) & = & a_{11}(\omega)x_1+a_{12}(\omega)x_2+
\widetilde{f}_1(\omega,x_1,x_2),\\
f_2(\omega,x_1,x_2) & = & a_{21}(\omega)x_1+a_{22}(\omega)x_2+
\widetilde{f}_2(\omega,x_1,x_2),
\end{eqnarray*}
where $\{a_{ij}(\omega)\}_{i,j=1,2}\in GL_2({\Bbb{Z}})$ and
$\widetilde{f}_1,\widetilde{f}_2:\Omega\times\tor^2\rightarrow\Bbb{R}$.
Then
\[u_1(\omega,x+1)=(a_{11}(\omega)-1)(x+1)+
\widetilde{f}_1(\omega,x+1,0)=u_1(\omega,x)+a_{11}(\omega)-1\]
and
\[u_1(\omega,x+c)=(a_{11}(\omega)-1)x-a_{12}(\omega)+
\widetilde{f}_1(\omega,x,-1)=u_1(\omega,x)-a_{12}(\omega).\]
Therefore
\[a_{11}(\omega)-1=\lim_{x\rightarrow+\infty}\frac{u_1(\omega,x)}{x}=
-a_{12}(\omega)/c\]
for $\mu$--a.e.\ $\omega\in\Omega$.
Since $c$ is irrational, we conclude that
$a_{11}(\omega)-1=a_{12}(\omega)=0$, hence that $u_1(\omega,\cdot)$ is $1$
and $c$ periodic, and finally that $u_1(\omega,\cdot)$ is a constant function
for $\mu$--a.e.\ $\omega\in\Omega$.
It is clear that the same conclusion can be drawn for $u_2$, which
completes the proof. $\Box$

\begin{lem}\label{rat}
If $c$ is rational, then there exist a group automorphism
$A:\tor^2\rightarrow\tor^2$, a $C^r$ random diffeomorphism of the circle
$F:\Omega\times\tor\rightarrow\tor$ and a $C^r$ random function
$\varphi:\Omega\times\tor\rightarrow\tor$ such that
\[A\circ f_{\omega}\circ A^{-1}(x_1,x_2)=(F_{\omega}x_1,
x_2+\varphi_{\omega}(x_1)).\] Moreover,
\begin{equation}\label{ht}
h_{T\omega}\circ
A^{-1}(F_{\omega}(x_1),x_2+\varphi_{\omega}(x_1))DF_{\omega}(x_1)=
h_{\omega}\circ A^{-1}(x_1,x_2)
\end{equation}
 for $\hat{\mu}$--a.e.\
$(\omega,x_1,x_2)\in\Omega\times\tor^2$, where
$\hat{\mu}:=(\mbox{\em Id}_{\Omega}\times A)\mu$ and
$h_{\omega}\circ A^{-1}:\tor^2\rightarrow{\Bbb{R}}$ depends only
on the first coordinate.
\end{lem}

\pf Let $p$ and $q$ be integer numbers such that $q>0$, $\gcd(p,q)=1$ and
$c=p/q$. Choose $a,b\in{\Bbb{Z}}$ with $ap-bq=1$. Consider the group
automorphism $A:\tor^2\rightarrow\tor^2$ associated to
the matrix
$A=\left[\begin{array}{rr}
q & -p\\
-b & a
\end{array}\right]$. Then $A^{-1}=\left[\begin{array}{rr}
a & p\\
b & q
\end{array}\right]$. Let us consider the $C^r$ RDS
$\hat{f}_{\omega}=A\circ f_{\omega}\circ A^{-1}$. Then $\hat{\mu}$ is a
$\hat{f}$--invariant measure and
\[D\hat{f}_{\omega}(x)=A\,Df_{\omega}(A^{-1}x)\,A^{-1}.\]
>From (\ref{piec}) we see that
\[\left[
\begin{array}{c}
p \\
q
\end{array}\right]
=Df_{\omega}x
\left[
\begin{array}{c}
p \\
q
\end{array}\right]
\]
for $P$--a.e.\ $\omega\in\Omega$ and all $x\in\tor^2$,
hence that
\[\left[
\begin{array}{c}
0 \\
1
\end{array}\right]
=D\hat{f}_{\omega}x
\left[
\begin{array}{c}
0 \\
1
\end{array}\right]
\]
for $P$--a.e.\ $\omega\in\Omega$ and all $x\in\tor^2$.
>From (\ref{szesc}) we see that
\[h_{\omega}(x)\left[
\begin{array}{cc}
q & -p
\end{array}\right]=
h_{T\omega}(f_{\omega}x)
\left[
\begin{array}{cc}
q & -p
\end{array}\right]
Df_{\omega}x
\]
for $\mu$--a.e.\ $(\omega,x)\in\Omega\times\tor^2$,
hence that
\[h_{\omega}\circ A^{-1}(x)\left[
\begin{array}{cc}
1 & 0
\end{array}\right]=
h_{T\omega}\circ A^{-1}(\hat{f}_{\omega}x)
\left[
\begin{array}{cc}
1 & 0
\end{array}\right]
D\hat{f}_{\omega}x
\]
for $\hat{\mu}$--a.e.\ $(\omega,x)\in\Omega\times\tor^2$.
It follows that
\[\frac{\partial(\hat{f}_{\omega})_1}{\partial x_2}=0,\;\;\;
\frac{\partial(\hat{f}_{\omega})_2}{\partial x_2}=1\]
for $P$--a.e.\ $\omega\in\Omega$  and
\[h_{T\omega}\circ
A^{-1}\circ\hat{f}_{\omega}\;\frac{\partial(\hat{f}_{\omega})_1}{\partial
x_1}=h_{T\omega}\circ A^{-1}\]
for $\hat{\mu}$--a.e.\ $(\omega,x)\in\Omega\times\tor^2$.
Therefore
\[\hat{f}_{\omega}(x_1,x_2)=(F_{\omega}x_1,
x_2+\varphi_{\omega}(x_1)),\]
where $F,\varphi:\Omega\times\tor\rightarrow\tor$ are $C^r$ random functions
and
\[h_{T\omega}\circ A^{-1}(F_{\omega}(x_1),x_2+\varphi_{\omega}(x_1))DF_{\omega}(x_1)=
h_{\omega}\circ A^{-1}(x_1,x_2)\] for  $\hat{\mu}$--a.e.\
$(\omega,x_1,x_2)\in\Omega\times\tor^2$. Since
$\hat{f}_{\omega}:\tor^2\rightarrow\tor^2$ is a
$C^r$--diffeomorphism, we
conclude that $F_{\omega}:\tor\rightarrow\tor$ is a
$C^r$--diffeomorphism for $P$--a.e.\ $\omega\in\Omega$. Since
\[\frac{1}{n^{\tau}}Df^n_{\omega}(x)\rightarrow h_{\omega}(x)
\left[
\begin{array}{c}
p/q \\ 1
\end{array}\right]
\left[
\begin{array}{cc}
1 & -p/q
\end{array}\right]
\]
for  $\mu$--a.e.\ $(\omega,x_1,x_2)\in\Omega\times\tor^2$, we have
\[\frac{1}{n^{\tau}}D\hat{f}^n_{\omega}(x)\rightarrow
h_{\omega}(A^{-1}x)/q^2 \left[\begin{array}{rr} 0 & 0\\ 1 & 0
\end{array}\right]
\]
for  $\hat{\mu}$--a.e.\ $(\omega,x_1,x_2)\in\Omega\times\tor^2$.
Let $\hat{h}_{\omega}:=h_{\omega}\circ A^{-1}$. Then
\[\frac{1}{n^{\tau}}\sum_{k=0}^{n-1}D\varphi_{T^k\omega}(F_{\omega}^{k}(x_1))
DF_{\omega}^k(x_1)\rightarrow\hat{h}_{\omega}(x_1,x_2)/q^2\] for
$\hat{\mu}$--a.e.\ $(\omega,x_1,x_2)\in\Omega\times\tor^2$. It
follows that $\hat{h}_{\omega}$ depends only on the first
coordinate. $\Box$

\vspace{2ex}

{\bf Proof of Theorem~\ref{rtg}.} By Lemmas~\ref{irrat} and
\ref{rat}, to prove the first claim of the theorem we only need to
show that $\tau \geq1$. Suppose that $\tau <1$. Let
$\nu:=(\mbox{Id}_{\Omega}\times\pi)\hat{\mu}$, where
$\pi:\tor^2\rightarrow\tor$
is the projection on the first coordinate. Then $\nu$ is a
$F$--invariant ergodic measure of full support on $\Omega\times\tor$. By
Lemma~\ref{rat},
\[\hat{h}_{T^k\omega}(F^k_{\omega}(x))DF^k_{\omega}(x)=
\hat{h}_{\omega}(x)\] and
\begin{equation}\label{doh}
\frac{1}{n^{\tau}}\sum_{k=0}^{n-1}D\varphi_{T^k\omega}(F_{\omega}^{k}(x))
DF_{\omega}^k(x)\rightarrow\hat{h}_{\omega}(x)/q^2
\end{equation}
for ${\nu}$--a.e.\ $(\omega,x)\in\Omega\times\tor$. Therefore
\begin{equation}\label{doje}
\frac{1}{n^{\tau}}\sum_{k=0}^{n-1}D\varphi_{T^k\omega}(F_{\omega}^{k}(x))/
\hat{h}_{T^k\omega}(F^k_{\omega}(x)) \rightarrow 1/q^2
\end{equation}
and consequently
\[\frac{1}{n}\sum_{k=0}^{n-1}D\varphi_{T^k\omega}(F_{\omega}^{k}(x))/
\hat{h}_{T^k\omega}(F^k_{\omega}(x)) \rightarrow 0\] for
${\nu}$--a.e.\ $(\omega,x)\in\Omega\times\tor$. It follows that
the measurable cocycle
$D\varphi/\hat{h}:\Omega\times\tor\rightarrow{\Bbb{R}}$ over the
skew product $T_F$ is recurrent (see \cite{Sch}). Consequently, for
${\nu}$--a.e.\
$(\omega,x)\in\Omega\times\tor$ there exists an increasing
sequence of natural numbers $\{n_i\}_{i\in{\Bbb{N}}}$ such that
\[|\sum_{k=0}^{n_i-1}D\varphi_{T^k\omega}(F_{\omega}^{k}(x))/
\hat{h}_{T^k\omega}(F^k_{\omega}(x))|\leq 1.\] It follows that
\[\frac{1}{n_i^{\tau}}\sum_{k=0}^{n_i-1}D\varphi_{T^k\omega}(F_{\omega}^{k}(x))/
\hat{h}_{T^k\omega}(F^k_{\omega}(x)) \rightarrow 0,\] contrary to
(\ref{doje}).

Now suppose that $f$ preserves orientation. Let us decompose
$\nu_{\omega}=\nu_{\omega}^d+\nu_{\omega}^c$, where
$\nu_{\omega}^d$ is the discrete and $\nu_{\omega}^c$ is the
continuous part of the measure $\nu_{\omega}$. As this
decomposition is measurable we can consider the measures
$\nu^d=\int_{\Omega}\nu_{\omega}^ddP\omega$ and
$\nu^c=\int_{\Omega}\nu_{\omega}^cdP\omega$ on $\Omega\times\tor$.
It is easy to check that $\nu^d$ and $\nu^c$ are $F$--invariant.
By the ergodicity of $\nu$, either $\nu=\nu^d$ or $\nu=\nu^c$.

We now show that $\nu=\nu^c$. Suppose, contrary to our claim, that
$\nu=\nu^d$. Let
$\Delta:\Omega\times\tor\rightarrow[0,1]$ denote the measurable function
given by
$\Delta(\omega,x)=\nu_{\omega}(\{x\})$. As $\nu$ is $F$--invariant we have
\[\Delta(T\omega,F_{\omega}x)=\nu_{T\omega}(\{F_{\omega}x\})=
F_{\omega}^{-1}\nu_{T\omega}(\{x\})=\nu_{\omega}(\{x\})=\Delta(\omega,x)\]
and consequently $\Delta$ is $T_F$--invariant. By the ergodicity of $T_F$,
the function $\Delta$ is $\nu$ constant. It follows that the measure
$\nu_{\omega}$ has only finite number of atoms for $P$--a.e.\
$\omega\in\Omega$, which contradicts the fact that $\nu$ has full support.

Define $\xi_{\omega}(x):=\int_0^xd\nu_{\omega}$ for all $x\in{\Bbb{R}}$.
Then  $\xi_{\omega}(x+1)=\xi_{\omega}(x)+1$, because
$\int_x^{x+1}d\nu_{\omega}=1$. Since $\nu_{\omega}$ is continuous and
$\nu$ has full support, the function
$\xi_{\omega}:{\Bbb{R}}\rightarrow{\Bbb{R}}$ is continuous and strictly
increasing. Therefore $\xi:\Omega\times\tor\rightarrow\tor$ is a random
homeomorphism. As $\nu$ is $F$--invariant and $F$ preserves orientation we
have
\[\xi_{T\omega}(F_{\omega}x)=\int_0^{F_{\omega}x}d\nu_{T\omega}=
\int_0^{F_{\omega}0}d\nu_{T\omega}+
\int_{F_{\omega}0}^{F_{\omega}x}dF_{\omega}\nu_{\omega}=
\alpha_{\omega}+\int_{0}^{x}d\nu_{\omega}=\xi_{\omega}(x)+\alpha_{\omega}\]
for $P$--a.e.\ $\omega\in\Omega$, where $\alpha_{\omega}=\int_0^{F_{\omega}0}d\nu_{T\omega}$.
$\Box$

\vspace{2ex}

{\bf Proof of Theorem~\ref{rtg1}} Suppose that $f$ has
$\tau$--polynomial $L^1$--growth of the derivative and $\mu$ is
equivalent to $P\otimes\lambda^{\otimes 2}$. Then $DF,D\varphi\in
L^1(\Omega\times\tor,\nu)$ and $\hat{\mu}$ is equivalent to
$P\otimes\lambda^{\otimes 2}$. Let $\theta\in
L^1(\Omega\times\tor^2,P\otimes\lambda^{\otimes 2})$ denote the
Radon--Nikodym derivative
of $\hat{\mu}$ with respect to $P\otimes\lambda^{\otimes 2}$. Then
\[\theta_{T\omega}(F_{\omega}(x_1),x_2+\varphi_{\omega}(x_1))DF_{\omega}(x_1)=
\theta_{\omega}(x_1,x_2)\] for  $P\otimes\lambda^{\otimes 2}$--a.e.\
$(\omega,x_1,x_2)\in\Omega\times\tor^2$. From (\ref{ht}) there exists a
non--zero constant $C$ such that
$\theta_{\omega}(x_1,x_2)=C\hat{h}_{\omega}(x_1)$ for
$P\otimes\lambda^{\otimes 2}$--a.e.\ $(\omega,x_1,x_2)\in\Omega\times\tor^2$.
Then the random homeomorphism $\xi_{\omega}:\tor\rightarrow\tor$ given by
$\xi_{\omega}(x):=\int_0^xd\nu_{\omega}=\int_0^x\theta_{\omega}(t)dt$ is a
Lipschitz random diffeomorphism, because $\theta$ and $1/\theta$ are
bounded. It follows that $f$ is Lipschitz
conjugate to the random skew product
\[(T_{\alpha,\psi})_{\omega}(x_1,x_2)=(x_1+\alpha_{\omega},
x_2+\psi_{\omega}(x_1)),\] where
$\psi_{\omega}:=\varphi_{\omega}\circ\xi^{-1}_{\omega}$. From (\ref{doh})
we conclude that $T_{\alpha,\psi}$ has $\tau$--polynomial
$L^1$--growth of the derivative and
\[\frac{1}{n^{\tau}}\sum_{k=0}^{n-1}D\psi(T_{\alpha}(\omega,x))
\rightarrow\hat{h}_{\omega}\circ\xi_{\omega}^{-1}(x)D\xi_{\omega}^{-1}(x)/q^2\]
in $L^1(\Omega\times\tor,P\otimes\lambda)$. It follows that
\[n^{1-\tau}\int_{\Omega}d(\psi_{\omega})dP\omega\rightarrow
1/Cq^2,\]
and finally that $\tau=1$ and $\int_{\Omega}d(\psi_{\omega})dP\omega\neq 0$.
$\Box$

\section{Area--preserving diffeomorphisms of the $3$--torus}\indent\label{3d}

In this section we give a classification of area--preserving
ergodic diffeomorphisms of polynomial uniform growth of the derivative
on the $3$--torus, i.e.\ diffeomorphisms of  polynomial growth
of the derivative for which the sequence
$\{n^{-\tau}Df^n\}_{n\in{\Bbb{N}}}$ converges uniformly.
We first give a sequence of essential examples of such diffeomorphisms.
We will consider 2--steps skew products
$T_{\alpha,\beta,\gamma}:\tor^3\rightarrow\tor^3$ given by
\[T_{\alpha,\beta,\gamma}(x_1,x_2,x_3)=
(x_1+\alpha,x_2+\beta(x_1),x_3+\gamma(x_1,x_2)),\]
where $\alpha$ is irrational and $\beta:\tor\rightarrow\tor$,
$\gamma:\tor^2\rightarrow\tor$ are
of class $C^1$. We will denote by $d_i(\gamma)$ the topological degree of
$\gamma$ with respect to the $i$-th coordinate for $i=1,2$.
We will use the symbol $h_{x_1}$ to denote the partial derivative
$\partial h/\partial x_i$.

\vspace{2ex}

{\bf Example 1.} Assume that $\beta$ is a constant function,
$\alpha,\beta,1$ are rationally independent and
$(d_1(\gamma),d_2(\gamma))\neq 0$. Then
\[\frac{1}{n}DT_{\alpha,\beta,\gamma}^n\rightarrow
\left[\begin{array}{ccc}
0 & 0 & 0 \\
0 & 0 & 0 \\
d_1(\gamma) & d_2(\gamma) & 0
\end{array}\right]\neq 0
\]
uniformly and  $T_{\alpha,\beta,\gamma}$
is ergodic, by Lemma~\ref{skew}.

\vspace{2ex}

{\bf Example 2.} Assume that $\gamma$ depends only on the first coordinate,
$d(\beta)=0$, $d(\gamma)\neq 0$ and the skew product
\[\tor^2\ni(x_1,x_2)\longmapsto(x_1+\alpha,x_2+\beta(x_1))\in\tor^2\]
is ergodic. Then
\[\frac{1}{n}DT_{\alpha,\beta,\gamma}^n\rightarrow
\left[\begin{array}{ccc}
0 & 0 & 0 \\
0 & 0 & 0 \\
d_1(\gamma) & 0 & 0
\end{array}\right]\neq 0
\]
uniformly. To prove the ergodicity $T_{\alpha,\beta,\gamma}$,
it suffices to show the ergodicity of the skew product
\[\tor^2\ni(x_1,x_2)\longmapsto(x_1+\alpha,x_2+m_1\beta(x_1)+m_2\gamma(x_1))\in\tor^2\]
for every $(m_1,m_2)\in{\Bbb{Z}}^2\setminus\{(0,0)\}$. If $m_2=0$, then
the ergodicity of the skew product follows from assumption. Otherwise, the
claim follows from Lemma~\ref{skew}.

\vspace{2ex}

{\bf Example 3.} Assume that
$d(\beta)\neq 0$  and $d_2(\gamma)\neq 0$.
By Lemma~\ref{skew}, $T_{\alpha,\beta,\gamma}$ is ergodic. We now show
that $T_{\alpha,\beta,\gamma}$ has square
uniform growth of the derivative. Precisely, we show that
\begin{equation}\label{kwa}
\frac{1}{n^2}DT_{\alpha,\beta,\gamma}^n\rightarrow
\left[\begin{array}{ccc}
0 & 0 & 0 \\
0 & 0 & 0 \\
d(\beta)d_2(\gamma)/2 & 0 & 0
\end{array}\right]\neq 0
\end{equation}
uniformly.
Let
\[T_{\alpha}(x)=x+\alpha\mbox{ and }T_{\alpha,\beta}(x_1,x_2)=(x_1+\alpha,x_2+\beta(x_1)).\]
Note that
\[T_{\alpha,\beta,\gamma}^n(x_1,x_2,x_3)=(x_1+n\alpha,x_2+\beta^{(n)}(x_1),x_3+\gamma^{(n)}(x_1,x_2)).\]
Therefore
\begin{eqnarray*}
\lefteqn{\frac{1}{n^2}DT_{\alpha,\beta,\gamma}^n(x_1,x_2,x_3)}\\ & &=
\left[\begin{array}{ccc}
n^{-2} & 0 & 0\\
n^{-2}(D\beta)^{(n)}(x_1) & n^{-2} & 0\\
n^{-2}(\gamma_{x_1}^{(n)}(x_1,x_2)+
\sum_{k=0}^{n-1}\gamma_{x_2}(T^k_{\alpha,\beta}(x_1,x_2))
(D\beta)^{(k)}(x_1)) &
n^{-2}(\gamma_{x_2})^{(n)}(x_1,x_2)
& n^{-2}
\end{array}\right]
\end{eqnarray*}
To obtain (\ref{kwa}), it only
remains to verify
\[\frac{1}{n^2}\sum_{k=0}^{n-1}\gamma_{x_2}
(T^k_{\alpha,\beta}(x_1,x_2))(D\beta)^{(k)}(x_1)
\rightarrow d(\beta)d_2(\gamma)/2\neq 0\]
uniformly. First observe that
\[\frac{1}{n^2}\sum_{k=0}^{n-1}(n-k)\gamma_{x_2}
(T^k_{\alpha,\beta}(x_1,x_2))=
\frac{1}{n^2}\sum_{k=1}^{n}(\gamma_{x_2})^{(k)}(x_1,x_2)
\rightarrow d_2(\gamma)/2\]
uniformly. Indeed, set $a_n:=\sum_{k=1}^{n}
(\gamma_{x_2})^{(k)}$ and $b_n:=n^2$. Then
\[\frac{a_{n+1}-a_n}{b_{n+1}-b_n}=\frac{\sum_{k=0}^{n-1}
\gamma_{x_2}\circ T^{k}_{\alpha,\beta}}{2n+1}
\rightarrow \frac{1}{2}\int_{\tor^2}
\gamma_{x_2}(x_1,x_2)dx_1dx_2=d_2(\gamma)/2\]
uniformly. It follows that $a_n/b_n\rightarrow d(\beta)d_2(\gamma)/2$
uniformly. Hence
\[\frac{1}{n^2}\sum_{k=0}^{n-1}k\gamma_{x_2}
\circ T^k_{\alpha,\beta}=
\frac{1}{n}\sum_{k=0}^{n-1}\gamma_{x_2}
\circ T^k_{\alpha,\beta}-
\frac{1}{n^2}\sum_{k=0}^{n-1}(n-k)\gamma_{x_2}
\circ T^k_{\alpha,\beta}\rightarrow d_2(\gamma)/2\]
uniformly. Then
\begin{eqnarray*}
\lefteqn{\lim_{n\rightarrow\infty}
\|\frac{1}{n^2}\sum_{k=0}^{n-1}\gamma_{x_2}
\circ T^k_{\alpha,\beta}\cdot D\beta^{(k)}-
d(\beta)d_2(\gamma)/2\|_{\sup}}\\
 & = &
\lim_{n\rightarrow\infty}\|\frac{1}{n^2}\sum_{k=0}^{n-1}\gamma_{x_2}
\circ T^k_{\alpha,\beta}\cdot D\beta^{(k)}-
\frac{1}{n^2}\sum_{k=0}^{n-1}kd(\beta)\gamma_{x_2}
\circ T^k_{\alpha,\beta}\|_{\sup}\\
 & \leq &
\lim_{n\rightarrow\infty}
\frac{\|\gamma\|_{C^1}}{n}
\sum_{k=0}^{n-1}\|\frac{1}{k}\sum_{j=0}^{k-1}D\beta\circ T_{\alpha}^{j}-d(\beta)\|_{\sup}
\\ & \leq &
\lim_{n\rightarrow\infty}
\frac{\|\gamma\|_{C^1}}{n}
\sum_{k=0}^{n-1}\|\frac{1}{k}\sum_{j=0}^{k-1}D\beta\circ T_{\alpha}^{j}-
\int_{\tor}D\beta(x)dx\|_{\sup}= 0.
\end{eqnarray*}

\vspace{2ex}

{\bf Example 4.} Assume that
\[f(x_1,x_2,x_3)=
(x_1+\alpha,-x_2+\beta(x_1),x_3+\gamma(x_1)),\]
where $\alpha$ is irrational and $\beta,\gamma:\tor\rightarrow\tor$
are of class $C^1$. Suppose that $d(\gamma)\neq 0$ and the factor map
\[\hat{f}:\tor^2\rightarrow\tor^2,\;\;\hat{f}(x_1,x_2)=(x_1+\alpha,-x_2+\beta(x_1))\]
is ergodic. It is easy to see that $\hat{f}$ is ergodic iff
$\hat{f}^2(x_1,x_2)=(x_1+2\alpha,x_2-\beta(x_1)+\beta(x_1+\alpha))$ is.
Since $d(-\beta+\beta\circ T_{\alpha})=0$ and $d(\gamma+\gamma\circ
T_{\alpha})\neq 0$, we conclude the diffeomorphism
\[f^2(x_1,x_2,x_3)=
(x_1+2\alpha,x_2-\beta(x_1)+\beta(x_1+\alpha),x_3+\gamma(x_1)+\gamma(x_1+\alpha))\]
is ergodic, by Example 2. Consequently, $f$ is ergodic. Moreover
\[\frac{1}{n}Df^n\rightarrow
\left[\begin{array}{ccc}
0 & 0 & 0 \\
0 & 0 & 0 \\
d(\gamma) & 0 & 0
\end{array}\right]\neq 0
\]
uniformly.

\vspace{2ex}

{\bf Example 5.} Assume that an area--preserving $C^2$--diffeomorphism $f:\tor^3\rightarrow\tor^3$ has
$\tau$--polynomial uniform growth of the derivative and
\[\frac{1}{n^{\tau}}Df^n\rightarrow
\left[\begin{array}{ccc}
0 & 0 & 0 \\
0 & 0 & 0 \\
h_1 & h_2 & 0
\end{array}\right]\neq 0,
\]
where $h=(h_1,h_2):\tor^3\rightarrow{\Bbb{R}}^2$ is of class $C^1$. Let
$\varphi:\tor^2\rightarrow\tor^2$ be an area--preserving
$C^2$--diffeomorphism. Let $\hat{f}:\tor^3\rightarrow\tor^3$ denote the
area--preserving diffeomorphism that is conjugate to $f$ via the
diffeomorphism $\varphi\times\mbox{Id}:\tor^3\rightarrow\tor^3$, i.e.
\[\hat{f}=(\varphi\times\mbox{Id})^{-1}\circ f\circ (\varphi\times\mbox{Id}).\]
Then
\[\frac{1}{n^{\tau}}D\hat{f}^n(x_1,x_2,x_3)\rightarrow
\left[\begin{array}{cc}
\begin{array}{cc}
 0 & 0
\end{array}
& 0 \\
\begin{array}{cc}
 0 & 0
\end{array}
 & 0 \\
Dh(\varphi(x_1,x_2),x_3)D\varphi(x_1,x_2) & 0
\end{array}\right]\neq 0
\]
uniformly. It follows that $\hat{f}$ has
$\tau$--polynomial uniform growth of the derivative and the limit function
$\lim_{n\rightarrow\infty}n^{-\tau}D\hat{f}^n$ is of class $C^1$.

\vspace{2ex}

In this section we prove the following theorem.

\begin{theo}\label{twgl}
Let $f:\tor^3\rightarrow\tor^3$ be an area--preserving ergodic
$C^2$--diffeomorphism with $\tau$--polynomial uniform growth of the
derivative ($\tau>0$). Suppose that the limit function
$\lim_{n\rightarrow\infty}n^{-\tau}Df^n$ is of class $C^1$. Then either
$\tau=1$ or $\tau=2$ and $f$ is $C^2$--conjugate to a diffeomorphism of
the form
\[\tor^3\ni(x_1,x_2,x_3)\longmapsto (x_1+\alpha,\ep
x_2+\beta(x_1),x_3+\gamma(x_1,x_2))\in\tor^3,\]
where $\ep=\det Df=\pm 1$.
\end{theo}

As in the previous section,
the proof of the main theorem is divided into a sequence of lemmas.
Suppose that $f:\tor^3\rightarrow\tor^3$ is an area--preserving ergodic
diffeomorphism with $\tau$--polynomial growth of the derivative. Let
$g:\tor^3\rightarrow M_3({\Bbb{R}})$ denote the limit of the sequence
$\{n^{-\tau}Df^n\}_{n\in{\Bbb{N}}}$. By Lemma~\ref{lemg},
$g(\bar{x})g(\bar{y})=0$ and $g(\bar{x})^2=0$ for all
$\bar{x},\bar{y}\in\tor^3$.

\begin{lem}\label{formag}
Suppose that $A,B$ are non--zero real $3\times 3$--matrixes such that
$A^2=B^2=AB=BA=0$. Then there exist three non--zero vectors (real $1\times
3$--matrixes) $\bar{a}$, $\bar{b}$, $\bar{c}$ such that
\begin{itemize}
\item $A=\bar{a}^T\,\bar{b}$ and $B=\bar{a}^T\,\bar{c}$ where
$\bar{b}\bar{a}^T=0$ and $\bar{c}\,\bar{a}^T=0$ or
\item $A=\bar{a}^T\,\bar{c}$ and $B=\bar{b}^T\,\bar{c}$  where
$\bar{c}\bar{a}^T=0$ and $\bar{c}\,\bar{b}^T=0$.
\end{itemize}
\end{lem}

\pf Suppose that $\lambda\in{\Bbb{C}}$ is an eigenvalues of $A$ for an
eigenvector $\bar{x}\in{\Bbb{C}}^3$. Then $\lambda^2\bar{x}=A^2\bar{x}=0$
and consequently $\lambda=0$. It follows that Jordan canonical form of
$A$ equals either
\[
\left[\begin{array}{ccc}
0 & 0 & 0 \\
1 & 0 & 0 \\
0 & 0 & 0
\end{array}\right] \mbox{ or }
\left[\begin{array}{ccc}
0 & 0 & 0 \\
1 & 0 & 0 \\
0 & 1 & 0
\end{array}\right].\]
But the second case can not occur since
\[\left[\begin{array}{ccc}
0 & 0 & 0 \\
1 & 0 & 0 \\
0 & 1 & 0
\end{array}\right]^2=
\left[\begin{array}{ccc}
0 & 0 & 0 \\
0 & 0 & 0 \\
1 & 0 & 0
\end{array}\right]\neq 0.\]
It follows that there exists $C\in GL_3({\Bbb{R}})$ such that
\[A=C
\left[\begin{array}{ccc}
0 & 0 & 0 \\
1 & 0 & 0 \\
0 & 0 & 0
\end{array}\right]C^{-1}
=\left[\begin{array}{c}
c_{12} \\
c_{22} \\
c_{32}
\end{array}\right]
\left[\begin{array}{ccc}
c_{11}^{-1} & c_{12}^{-1} & c_{13}^{-1}
\end{array}\right].
\]
Therefore we can find non--zero real $1\times
3$--matrixes $\bar{a}_1$, $\bar{a}_2$ such that
$A=\bar{a}_1^T\,\bar{a}_2$. As $A^2=0$ we have $\bar{a}_1\perp\bar{a}_2$.
Similarly,
we can find non--zero real $1\times
3$--matrixes $\bar{b}_1$,
$\bar{b}_2$ such that
$B=\bar{b}_1^T\,\bar{b}_2$ and $\bar{b}_1\perp\bar{b}_2$. Let
$\bar{o}\in{\Bbb{R}}^3$ be a non--zero vector that is orthogonal to both
$\bar{a}_1$ and $\bar{a}_2$. As $AB=BA=0$ we have $\bar{a}_1\perp\bar{b}_2$ and
$\bar{a}_2\perp\bar{b}_1$. It follows that there exists a real matrix
$[d_{ij}]_{i,j=1,2}$ such that
\begin{eqnarray*}
\bar{b}_1 & = & d_{11}\bar{a}_1+d_{12}\bar{o}\\
\bar{b}_2 & = & d_{21}\bar{a}_2+d_{22}\bar{o}.
\end{eqnarray*}
Then
\[0=\langle\bar{b}_1,\bar{b}_2\rangle=d_{12}d_{22}\|\bar{o}\|^2.\]
If $d_{12}=0$, then $d_{11}\neq 0$ and
we can put $\bar{a}:=\bar{a}_1$, $\bar{b}:=\bar{a}_2$,
$\bar{c}:=d_{11}\bar{b}_2$. Then  $\bar{a}^T\,\bar{b}=A$ and
$\bar{a}^T\,\bar{c}=B$.
If $d_{22}=0$, then $d_{21}\neq 0$ and
we can put $\bar{a}:=\bar{a}_1/d_{21}$, $\bar{b}:=\bar{b}_1$,
$\bar{c}:=\bar{b}_2$. Then $\bar{a}^T\,\bar{c}=A$ and
$\bar{b}^T\,\bar{c}=B$, which completes the proof. $\Box$

\vspace{2ex}

By the above lemma, there exists $\bar{c}\in{\Bbb{R}}^3$ such that
for any two linearly independent vectors $\bar{a},\bar{b}\in{\Bbb{R}}^3$
orthogonal to $\bar{c}$ there exist $C^1$--functions
$h_1,h_2:\tor^3\rightarrow {\Bbb{R}}$ such that $g(\bar{x})$ equals
\[\bar{c}^T(h_1(\bar{x})\bar{a}+h_2(\bar{x})\bar{b})\mbox{ or }
(h_1(\bar{x})\bar{a}+h_2(\bar{x})\bar{b})^T\bar{c}\]
for all $\bar{x}\in\tor^3$.
First we deal with the special case where the limit function $g$ is
constant.

\begin{lem}\label{staleg}
Let $f:\tor^3\rightarrow\tor^3$ be an area--preserving ergodic
$C^1$--diffeomorphism with $\tau$--polynomial uniform growth of the
derivative ($\tau >0$). Suppose that the limit function
$g=\lim_{n\rightarrow\infty}n^{-\tau}Df^n$ is constant. Then either
$\tau=1$ or $\tau=2$ and $f$ is algebraically conjugate to a
diffeomorphism of
the form
\[
\tor^3\ni(x_1,x_2,x_3)\longmapsto (x_1+\alpha,\ep
x_2+\beta(x_1),x_3+\gamma(x_1,x_2))\in\tor^3,\]
where $\ep=\det Df=\pm 1$.
\end{lem}

Before we will pass to the proof we give notation.
Let $A\in GL_3({\Bbb{R}})$.
Let us denote by $\tor^3_A$ the quotient group ${\Bbb R}^3/{\Bbb{Z}}^3A^T$,
which is a model of the $3$--torus as well. Then the map
\[A:\tor^3\rightarrow\tor^3_A,\;\;\;\;A\bar{x}=\bar{x}A^T\]
establishes a smooth isomorphism between $\tor^3$ and $\tor^3_A$.
Suppose that $\xi:\tor^3_A\rightarrow\tor^3_A$ is a diffeomorphism.
Then $A^{-1}\circ \xi\circ A$ is a diffeomorphism of the torus $\tor^3$.
Let $N\in GL_3({\Bbb{Z}})$ be its linear part. Then
\[\xi(\bar{x}+\bar{m}A^T)=\xi(\bar{x})+\bar{m}N^TA^T\]
for all $\bar{m}\in{\Bbb{Z}}^3$. Moreover, we can decompose
\[\xi(\bar{x})=\bar{x}(ANA^{-1})^T+\tilde{\xi}(\bar{x})\]
where $ANA^{-1}$ is the linear and $\tilde{\xi}$ is the periodic part of
$\xi$, i.e.\
\[\tilde{\xi}(\bar{x}+\bar{m}A^T)=\tilde{\xi}(\bar{x})\]
for all $\bar{m}\in{\Bbb{Z}}^3$.

Suppose that ${f}:\tor^3\rightarrow\tor^3$ is a smooth  diffeomorphism
with $\tau$--polynomial uniform growth of the derivative and $g:\tor^3\rightarrow
M_3({\Bbb{R}})$ is the limit of the sequence
$\{n^{-\tau}Df^n\}_{n\in{\Bbb{N}}}$.
Let us consider the diffeomorphism $\hat{f}:\tor^3_A\rightarrow\tor^3_A$
given by $\hat{f}:=A\circ f\circ A^{-1}$.
Then
\begin{equation}\label{betahat}
\frac{1}{n^{\tau}}D\hat{f}^n(\bar{x})=
\frac{1}{n^{\tau}}ADf^n(A^{-1}\bar{x})A^{-1}\rightarrow
Ag(A^{-1}\bar{x})A^{-1}
\end{equation}
uniformly on $\tor^3_A$. Let us denote by
$\hat{g}:\tor^3_A\rightarrow M_3({\Bbb{R}})$ the function
$\hat{g}(\bar{x}):=Ag(A^{-1}\bar{x})A^{-1}$.
Lemma~\ref{lemg} now gives
\begin{equation}\label{kapel0}
g(\bar{x})=g(f\bar{x})Df(\bar{x})\;\;\;\;\mbox{ and }\;\;\;\;
g(\bar{y})=Df(\bar{x})\,g(\bar{y})
\end{equation}
for all $\bar{x},\bar{y}\in\tor^3$, and consequently
\begin{eqnarray}\label{kapel}
\hat{g}(\bar{x})=\hat{g}(\hat{f}\bar{x})D\hat{f}(\bar{x})
 &\mbox{ and } &
\hat{g}(\bar{y})=D\hat{f}(\bar{x})\,\hat{g}(\bar{y})
\end{eqnarray}
for all $\bar{x},\bar{y}\in\tor^3_A$.

Throughout this paper we will denote by $G(\bar{c})$ the subgroup of all
$\bar{m}\in{\Bbb{Z}}^3$  such that $\bar{m}\perp\bar{c}$. Of course, the
rank of $G(\bar{c})$ can be equal $0$, $1$ or $2$, whenever $c\neq 0$.
Further useful properties of the group $G(\bar{c})$ the reader can find in
Appendix~\ref{app}.

\vspace{2ex}

Suppose that $f:\tor^3\rightarrow\tor^3$ is an area--preserving ergodic
$C^1$--diffeomorphism with $\tau$--polynomial uniform growth of the
derivative and the limit function $g$ is constant.
By Lemma~\ref{formag}, there exist mutually orthogonal vectors
$\bar{a},\bar{c}\in{\Bbb{R}}^3$ such that $g=\bar{c}^T\bar{a}$.
We now prove that following

\begin{lem}\label{pomoc}
Let $f:\tor^3\rightarrow\tor^3$ be an area--preserving
$C^1$--diffeomorphism. Suppose that $f$ preserves orientation,  has
$\tau$--polynomial uniform growth of the
derivative and the limit function
$g=\lim_{n\rightarrow\infty}n^{-\tau}Df^n$ equals $\bar{c}^T\bar{a}$. Then
rank $G(\bar{a})=2$. Moreover, $\tau$ equals either $1$ or $2$.
\end{lem}

\pf
Let $\bar{b}\in{\Bbb{R}}^3$ be a vector orthogonal to both $\bar{a}$ and
$\bar{c}$ such that $\det(A)=1$, where
\[
A=\left[\begin{array}{c}
\bar{a} \\
\bar{b} \\
\bar{c}
\end{array}\right].\]
Consider $\hat{f}:\tor^3_A\rightarrow\tor^3_A$
given by $\hat{f}:=A\circ f\circ A^{-1}$. Then
\[\hat{g}=
A\bar{c}^T\bar{a}
A^{-1}
=\left[\begin{array}{c}
0 \\
0 \\
1
\end{array}\right]
\left[\begin{array}{ccc}
1 &
0 &
0
\end{array}\right].\]
>From (\ref{kapel}) we obtain
\begin{eqnarray*}
\left[\begin{array}{ccc}
1 &
0 &
0
\end{array}\right]=
\left[\begin{array}{ccc}
1 &
0 &
0
\end{array}\right]
D\hat{f}
 &\mbox{ and } &
\left[\begin{array}{c}
0 \\
0 \\
1
\end{array}\right]
=D\hat{f}
\left[\begin{array}{c}
0 \\
0 \\
1
\end{array}\right].
\end{eqnarray*}
Consequently,
\begin{eqnarray*}
\frac{\partial}{\partial x_1}\hat{f}_1(\bar{x})=1, &
\displaystyle \frac{\partial}{\partial x_2}\hat{f}_1(\bar{x})=0, &
\frac{\partial}{\partial x_3}\hat{f}_1(\bar{x})=0,\\
\frac{\partial}{\partial x_3}\hat{f}_1(\bar{x})=0, &
\displaystyle \frac{\partial}{\partial x_3}\hat{f}_2(\bar{x})=0, &
\frac{\partial}{\partial x_3}\hat{f}_3(\bar{x})=1
\end{eqnarray*}
for all $\bar{x}\in\tor^3_A$. It follows that
\[\hat{f}(x_1,x_2,x_3)=(x_1+\alpha, x_2+\beta(x_1),x_3+
\gamma(x_1,x_2)),\]
where $\beta:{\Bbb{R}}\rightarrow{\Bbb{R}}$,
$\gamma:{\Bbb{R}}^2\rightarrow{\Bbb{R}}$ are $C^1$--function. Let
$N\in GL_3({\Bbb{Z}})$ denote the linear part of $f$.
Then the linear part of $\hat{f}$ equals
\[ANA^{-1}=
\left[\begin{array}{ccc}
1 & 0 & 0 \\
K_{21} & 1 & 0 \\
K_{31} & K_{32} & 1
\end{array}\right].
\]
It follows that
\begin{eqnarray}
\bar{a}N & = & \bar{a} \label{wla}\\
\bar{b}N & = & K_{21}\bar{a}+\bar{b} \label{wla1} \\
\bar{c}N & = & K_{31}\bar{a}+K_{32}\bar{b}+\bar{c}.\label{wla2}
\end{eqnarray}
Let $\tilde{f}:\tor^3\rightarrow{\Bbb{R}}^3$ denote the periodic part of
$f$, i.e.\
\[f(\bar{x})=\bar{x}N^T+\tilde{f}(\bar{x}).\]
Then
\[f^n(\bar{x})=\bar{x}(N^n)^T+\sum_{k=0}^{n-1}\tilde{f}(f^k\bar{x})(N^{n-1-k})^T.\]
Since $\int_{\tor^3}D(\tilde{f}\circ f^k)(\bar{x})d\bar{x}=0$ for all
natural $k$, we see that
\begin{equation}\label{cud}
\frac{1}{n^{\tau}}N^n=\frac{1}{n^{\tau}}\int_{\tor^3}
Df^n(\bar{x})d\bar{x}\rightarrow g.
\end{equation}
It follows that
\begin{equation}\label{cud1}
\frac{1}{n^{\tau}}\left[\begin{array}{ccc}
1 & 0 & 0 \\
K_{21} & 1 & 0 \\
K_{31} & K_{32} & 1
\end{array}\right]^n
\rightarrow \hat{g}
=\left[\begin{array}{ccc}
0 & 0 & 0\\
0 & 0 & 0\\
1 & 0 & 0
\end{array}\right]
.
\end{equation}
Suppose, contrary to our claim,  that  rank $G(\bar{a})<2$.

First suppose that rank $G(\bar{a})=0$. From (\ref{wla}) we conclude that
$N=$ Id, and finally that $n^{-\tau}N^n$ tends to zero, contrary to
(\ref{cud}).

Now suppose that rank $G(\bar{a})=1$. Let $\bar{m}\in{\Bbb{Z}}^3$ be a
generator of $G(\bar{a})$. Then there exists a rational vector
$\bar{r}\in{\Bbb{Q}}^3$ such that $N-\mbox{Id}=\bar{m}^T\bar{r}$, by
(\ref{wla}).
>From (\ref{wla1}) we have
\[\bar{b}\bar{m}^T\bar{r}=\bar{b}(N-\mbox{Id})=K_{21}\bar{a}.\]
Suppose that $K_{21}\neq 0$. Then rank $G(\bar{a})=$ rank $G(\bar{r})=2$,
which contradicts our assumption. Consequently, $K_{21}=0$. It follows that
\[\left[\begin{array}{ccc}
1 & 0 & 0 \\
K_{21} & 1 & 0 \\
K_{31} & K_{32} & 1
\end{array}\right]^n=\left[\begin{array}{ccc}
1 & 0 & 0 \\
0 & 1 & 0 \\
nK_{31} & nK_{32} & 1
\end{array}\right].\]
>From (\ref{cud1}) it follows that $n=1$ and $K_{31}=1, K_{32}=0$.
Then
\[\bar{c}\bar{m}^T\bar{r}=\bar{c}(N-\mbox{Id})=\bar{a},\]
by (\ref{wla2}). It follows that rank $G(\bar{a})=$ rank $G(\bar{r})=2$,
which contradicts our assumption.

Finally, we have to prove that $\tau$ equals either $1$ or $2$.
>From (\ref{cud1}) we obtain
\[n^{1-\tau}K_{21}\rightarrow 0,\;\;\;\;
n^{1-\tau}K_{31}+\frac{1-1/n}{2}n^{2-\tau}K_{21}K_{32}\rightarrow 1,\;\;\;\;
n^{1-\tau}K_{32}\rightarrow 0.\]
If $K_{21}=0$, then $\tau=1$ and $K_{31}=1$. Otherwise, $\tau=2$ and
$K_{21}K_{32}=2$, which completes the proof. $\Box$

\vspace{2ex}

{\bf Proof of Lemma~\ref{staleg}.}
First note that $f^2$ preserves area and orientation and
$n^{-\tau}Df^{2n}$ uniformly tends to $2^{\tau}\bar{c}^T\bar{a}$. By
Lemma~\ref{pomoc}, rank $G(\bar{a})=2$. It follows that $\bar{a}=a\bar{m}\in
a{\Bbb{Z}}^3$, by Lemma~\ref{per} (see Appendix~\ref{app}).
Now choose $\bar{n},\bar{k}\in{\Bbb{Z}}^3$ such that the determinant of
\[A:=\left[
\begin{array}{c}
\bar{m}\\
\bar{n}\\
\bar{k}
\end{array}\right]\]
equals $1$. Let us consider the diffeomorphism
$\hat{f}:\tor^3\rightarrow\tor^3$
given by $\hat{f}:=A\circ f\circ A^{-1}$. Then
\[\hat{g}=AgA^{-1}=a
\left[\begin{array}{c}
0 \\
\bar{n}\bar{c}^T \\
\bar{k}\bar{c}^T
\end{array}\right]
\left[\begin{array}{ccc}
1 & 0 & 0
\end{array}\right].
\]
>From (\ref{kapel}) we have
\begin{eqnarray*}
\left[\begin{array}{ccc}
1 & 0 & 0
\end{array}\right]D\hat{f}(\bar{x})=
\left[\begin{array}{ccc}
1 & 0 & 0
\end{array}\right]
 & \mbox{ and } &
\left[\begin{array}{c}
0 \\
\bar{n}\bar{c}^T \\
\bar{k}\bar{c}^T
\end{array}\right]
=D\hat{f}(\bar{x})
\left[\begin{array}{c}
0 \\
\bar{n}\bar{c}^T \\
\bar{k}\bar{c}^T
\end{array}\right].
\end{eqnarray*}
It follows that
\[\hat{f}(x_1,x_2,x_3)=(x_1+\alpha,\varphi_{x_1}(x_2,x_3)),\]
where $\varphi:\tor\times\tor^2\rightarrow\tor^2$ is an area--preserving
random diffeomorphism over the rotation by an irrational number $\alpha$.
Then
\[\left[\begin{array}{c}
\bar{n}\bar{c}^T \\
\bar{k}\bar{c}^T
\end{array}\right]
=D\varphi
\left[\begin{array}{c}
\bar{n}\bar{c}^T \\
\bar{k}\bar{c}^T
\end{array}\right].
\]

Suppose that $\bar{n}\bar{c}^T$ and $\bar{k}\bar{c}^T$ are rationally
independent. By Lemma~\ref{irrat},
$\varphi_{x_1}(x_2,x_3)=(x_2+\beta(x_1),x_3+\gamma(x_1))$, where
$\beta, \gamma:\tor\rightarrow\tor$ are $C^1$--functions, which is our
assertion.

Otherwise, by Lemma~\ref{rat}, there exist
a group automorphism
$B:\tor^2\rightarrow\tor^2$ and $C^1$--functions
$\beta:\tor\rightarrow\tor$, $\gamma:\tor^2\rightarrow\tor$ such that
\[B\circ \varphi_{x_1}\circ B^{-1}(x_2,x_3)=(\ep x_2+\beta(x_1),
x_3+\gamma(x_1,x_2)),\]
where $\ep=\det Df$, which implies our assertion.
$\Box$

\vspace{2ex}

{\bf Proof of Theorem~\ref{twgl}} is divided into a few cases.

\vspace{2ex}

{\em Case} 1: Suppose that $g=\bar{c}^T(h_1\bar{a}+h_2\bar{b})$, where
$\bar{a}$ and $\bar{b}$ are orthogonal to $\bar{c}$ and the matrix
\[
A=\left[\begin{array}{c}
\bar{a} \\
\bar{b} \\
\bar{c}
\end{array}\right]\]
is nonsingular.
Let $\hat{f}:\tor^3_A\rightarrow\tor^3_A$ be given by $\hat{f}:=A\circ
f\circ A^{-1}$.
Then
\[\hat{g}=
A\bar{c}^T(\hat{h}_1\bar{a}+\hat{h}_2\bar{b})
A^{-1}
=\left[\begin{array}{c}
0 \\
0 \\
1
\end{array}\right]
\left[\begin{array}{ccc}
\hat{h}_1 &
\hat{h}_2 &
0
\end{array}\right],\]
where $\hat{h}_i(\bar{x}):=h_i(A^{-1}\bar{x})$ for $i=1,2$.
>From (\ref{kapel}) we obtain
\begin{eqnarray}\label{el}
\left[\begin{array}{ccc}
\hat{h}_1(\bar{x}) &
\hat{h}_2(\bar{x}) &
0
\end{array}\right] & = &
\left[\begin{array}{ccc}
\hat{h}_1(\hat{f}\bar{x}) &
\hat{h}_2(\hat{f}\bar{x}) &
0
\end{array}\right]
D\hat{f}(\bar{x}),\\ \nonumber
\left[\begin{array}{c}
0 \\
0 \\
1
\end{array}\right]
 & = & D\hat{f}(\bar{x})
\left[\begin{array}{c}
0 \\
0 \\
1
\end{array}\right]
\end{eqnarray}
for all $\bar{x}\in\tor^3_A$. Consequently
\[\frac{\partial}{\partial x_3}\hat{f}_1(\bar{x})=0,\;\;\;\;
\frac{\partial}{\partial x_3}\hat{f}_2(\bar{x})=0\;\;\;\;
\mbox{ and }\;\;\;
\frac{\partial}{\partial x_3}\hat{f}_3(\bar{x})=1
\]
for all $\bar{x}\in\tor^3_A$. It follows that
\[\hat{f}(x_1,x_2,x_3)=(F(x_1,x_2),x_3+
\gamma(x_1,x_2)),\]
where $\gamma:{\Bbb{R}}^2\rightarrow{\Bbb{R}}$ is a smooth function and
$F:{\Bbb{R}}^2\rightarrow{\Bbb{R}}^2$ is the diffeomorphism given by
$F(x_1,x_2)=(\hat{f}_1(x_1,x_2),\hat{f}_2(x_1,x_2))$.
Let $K$ denote the linear part of $\hat{f}$. Then $K=A\,N\,A^{-1}$, where
$N\in GL_3({\Bbb{Z}})$ is the linear part of
$f$. Then $\det K=\det N=\ep'=\pm 1$ and $K_{13}=0,K_{23}=0,K_{33}=1$.
Moreover, there exist
$(\bar{a}\bar{m}^T,\bar{b}\bar{m}^T)$--periodic for all
$\bar{m}\in{\Bbb{Z}}^3$ functions
$\tilde{F}:{\Bbb{R}}^2\rightarrow{\Bbb{R}}^2$,
$\tilde{\gamma}:{\Bbb{R}}^2\rightarrow{\Bbb{R}}$ such that
\begin{eqnarray*}
F(\bar{x}) & = & \tilde{F}(\bar{x})+\bar{x}K'^T,\\
{\gamma}(x_1,x_2) & = & \tilde{\gamma}(x_1,x_2)+ K_{31}x_1+ K_{32}x_2,
\end{eqnarray*}
where $K'=K|_{\{1,2\}\times\{1,2\}}\in GL_2({\Bbb{R}})$ and $\det K'=\ep'$.
>From (\ref{betahat}) we deduce that
\begin{eqnarray*}
\frac{1}{n^{\tau}}DF^n(x_1,x_2)
\rightarrow 0 & \mbox{ and } &
\frac{1}{n^{\tau}}\sum_{k=0}^{n-1}D(\gamma\circ F^k)(x_1,x_2)
\rightarrow
\left[\,\hat{h}_1(\bar{x}) \;
\hat{h}_2(\bar{x})\,\right]\end{eqnarray*}
uniformly on $\tor^3_A$. Therefore $\hat{h}_1,\hat{h}_2$ depend only on
the first two coordinates.  Let $H:{\Bbb{R}}^2\rightarrow{\Bbb{R}}^2$ be
given by $H(x_1,x_2)=\left[\,\hat{h}_1(x_1,x_2,0)\;
\hat{h}_2(x_1,x_2,0)\,\right]$. Then $H$ is
$(\bar{a}\bar{m}^T,\bar{b}\bar{m}^T)$--periodic for all
$\bar{m}\in{\Bbb{Z}}^3$.
>From (\ref{el}) we have
\begin{equation}\label{mot}
H(F\bar{x})DF(\bar{x})=H(\bar{x})
\end{equation}
for all $\bar{x}\in{\Bbb{R}}^2$.
Set $\chi_n:=n^{-\tau}\sum_{k=0}^{n-1}\gamma\circ F^k$. Since
$D\chi_n\rightarrow H$ uniformly on ${\Bbb{R}}^2$, we see that
\begin{eqnarray*}
\chi_n(x_1,x_2)-\chi_n(x_1,0)  & \rightarrow & \int_0^{x_2}H_2(x_1,t)dt,\\
\chi_n(x_1,x_2)-\chi_n(0,x_2) & \rightarrow & \int_0^{x_1}H_1(t,x_2)dt
\end{eqnarray*}
for all $(x_1,x_2)\in{\Bbb{R}}^2$.
Let $\xi:{\Bbb{R}}^2\rightarrow{\Bbb{R}}$ be defined by
\[\xi(x_1,x_2):=\lim_{n\rightarrow\infty}(\chi_n(x_1,x_2)-\chi_n(0,0))=
\left\{
\begin{array}{l}
\int_0^{x_1}H_1(t,x_2)dt+\int_0^{x_2}H_2(0,t)dt\\
\int_0^{x_2}H_2(x_1,t)dt+\int_0^{x_1}H_1(t,0)dt.
\end{array}\right.\]
Then $\partial \xi/\partial x_1=H_1$, $\partial \xi/\partial x_2=H_2$
and $\xi$ is of class $C^2$. From (\ref{mot}) it follows that there exists
$\alpha\in{\Bbb{R}}$ such that
\begin{equation}\label{pryk}
\xi(F\bar{x})=\xi(\bar{x})+\alpha.
\end{equation}
By Lemma~\ref{per} (see Appendix~\ref{app}), there exist a
$(\bar{a}\bar{m}^T,\bar{b}\bar{m}^T)$--periodic for all
$\bar{m}\in{\Bbb{Z}}^3$ function
$\tilde{\xi}:{\Bbb{R}}^2\rightarrow{\Bbb{R}}$ and $d_1,d_2\in{\Bbb{R}}$
such that $\xi(x_1,x_2)=\tilde{\xi}(x_1,x_2)+d_1x_1+d_2x_2$. Since $H\neq 0$,
it is easy to see that $(d_1,d_2)\neq (0,0)$. Moreover, from (\ref{pryk}) we
have
\begin{equation}\label{syf0}
\left[\,d_1\;d_2\,\right]\;K'=\left[\,d_1\;d_2\,\right]
\end{equation}
and
\begin{eqnarray}\label{syf1}
\lefteqn{\tilde{\xi}(\bar{x})+\alpha}\\
& = & \nonumber
\tilde{\xi}(\tilde{F}_1(\bar{x})+K_{11}x_1+K_{12}x_2,
\tilde{F}_2(\bar{x})+K_{21}x_1+K_{22}x_2)+d_1\tilde{F}_1(\bar{x})
+d_2\tilde{F}_2(\bar{x}).
\end{eqnarray}

{\em Case} 1a: Suppose that rank $G(\bar{c})=0$. By Lemma~\ref{per}, $D\hat{f}$ is
constant. It follows that $Df$ and $g$ are constant. Therefore
$g=\bar{c}^T\bar{a}$, where $\bar{a}$ is orthogonal to $\bar{c}$. From
(\ref{kapel0}) we obtain $\bar{c}^T=Df\bar{c}^T$. As $G(\bar{c})=\{0\}$ we have
$Df=$Id. Consequently, $f$ is a rotation on the $3$--torus, which
is impossible.

\vspace{2ex}

{\em Case } 1b: Suppose that rank $G(\bar{c})=1$. By Lemma~\ref{per}, there exist
$\bar{F}:\tor\rightarrow{\Bbb{R}}^2$, $\bar{\xi}:\tor\rightarrow{\Bbb{R}}$,
 $\bar{\gamma}:\tor\rightarrow{\Bbb{R}}$ and real numbers $l_1$, $l_2$
such that $\bar{m}=l_1\bar{a}+l_2\bar{b}$ generates  $G(\bar{c})$
and
\begin{eqnarray*}
\tilde{F}(x_1,x_2) & = & \bar{F}(l_1x_1+l_2x_2),\\
\tilde{\xi}(x_1,x_2) & = & \bar{\xi}(l_1x_1+l_2x_2),\\
\tilde{\gamma}(x_1,x_2) & = & \bar{\gamma}(l_1x_1+l_2x_2).
\end{eqnarray*}
>From (\ref{syf1}) we obtain
\begin{eqnarray*}
\lefteqn{\bar{\xi}(l_1x_1+l_2x_2)+\alpha}\\
 & = &
\bar{\xi}(l_1\bar{F}_1(l_1x_1+l_2x_2)+
l_2\bar{F}_2(l_1x_1+l_2x_2)+s_{1}x_1+s_{2}x_2)\\
 & & +d_1\bar{F}_1(l_1x_1+l_2x_2)
+d_2\bar{F}_2(l_1x_1+l_2x_2),
\end{eqnarray*}
where $\left[\,s_1\;s_2\,\right]=\left[\,l_1\;l_2\,\right]K'$.
If $(s_1,s_2)$ and $(l_1,l_2)$ are linearly independent, then $\bar{\xi}$
is constant. It follows that $H$ is constant and we can come back to
Lemma~\ref{staleg}.
Otherwise, there exists a real number $s$ such that
$(s_1,s_2)=s(l_1,l_2)$ and
\[\bar{\xi}(x)+\alpha =
\bar{\xi}(l_1\bar{F}_1(x)+
l_2\bar{F}_2(x)+sx)+d_1\bar{F}_1(x)+d_2\bar{F}_2(x)\]
for any real $x$.
As $f$ preserves area we have $\det DF(\bar{x})=\ep=\pm 1$. It follows
that
\begin{eqnarray*}
\ep & = & \det\left[
\begin{array}{cc}
l_1D\bar{F}_1(x)+K_{11} & l_2D\bar{F}_1(x)+K_{12}\\
l_1D\bar{F}_2(x)+K_{21} & l_2D\bar{F}_2(x)+K_{22}
\end{array}\right]\\
 & = &
(l_1K_{22}-l_2K_{21})D\bar{F}_1(x)+(-l_1K_{12}+l_2K_{11})D\bar{F}_2(x)+
\det K\\
 & = &
(l_1D\bar{F}_1(x)+l_2D\bar{F}_2(x))\det K/s+\det K
\end{eqnarray*}
for any real $x$. Since $\bar{F}_1, \bar{F}_2$ are $1$--periodic, we have
$l_1D\bar{F}_1(x)+l_2D\bar{F}_2(x)=0$ and $\det K=\ep$.
Therefore the function $l_1\bar{F}_1+l_2\bar{F}_2$ is constant.
Let us choose real numbers $r_1,r_2$ such that the determinant of the
matrix
\[L=\left[
\begin{array}{ccc}
l_1 & l_2 & 0\\
r_1 & r_2 & 0\\
0 & 0 & 1
\end{array}\right]\]
equals $1$. Now consider the diffeomorphism
$\check{f}:\tor^3_{LA}\rightarrow\tor^3_{LA}$ given by
$\check{f}=L\circ\hat{f}\circ L^{-1}$. Then
\[\check{f}(x_1,x_2,x_3)=(sx_1+\alpha,
\ep/sx_2+rx_1+r_1\bar{F}_1(x_1)+r_2\bar{F}_2(x_1),
x_3+\bar{\gamma}(x_1)+p_1x_1+p_2x_2).\]
As $\partial\check{f}^n_1/\partial x_1=s^n$ and
$\partial\check{f}^n_2/\partial x_2=(\ep/s)^n$ we obtain $s=\pm
1$, because $\check{f}$ has polynomial uniform growth of the derivative.
Moreover,
\[LA=\left[
\begin{array}{c}
\bar{m}\\
r_1\bar{a}+r_2\bar{b}\\
\bar{c}
\end{array}\right]\]
and $LA\circ f=\check{f}\circ LA$. Therefore
$f(\bar{x})\bar{m}^T=s\bar{x}\bar{m}^T+\alpha$. Observe that $s=1$.
Indeed, suppose, contrary to our claim, that $s=-1$. Consider the smooth
function $\kappa:\tor^3\rightarrow{\Bbb{C}}$ given by
$\kappa(\bar{x})=e^{2\pi i\bar{x}\bar{m}^T}$. Then $\kappa\circ
f^2=\kappa$. Since $\kappa$ is smooth, we conclude that is constant, by
the ergodicity of $f$. Consequently, $\bar{m}=0$, which is impossible.

Now choose $\bar{n},\bar{k}\in{\Bbb{Z}}^3$ such that the determinant of
\[A:=\left[
\begin{array}{c}
\bar{m}\\
\bar{n}\\
\bar{k}
\end{array}\right]\]
equals $1$. Let us consider the diffeomorphism
$\hat{f}:\tor^3\rightarrow\tor^3$
given by $\hat{f}:=A\circ f\circ A^{-1}$.
>From (\ref{kapel}) we have
\[\left[\begin{array}{c}
0 \\
\bar{n}\bar{c}^T \\
\bar{k}\bar{c}^T
\end{array}\right]
=D\hat{f}(\bar{x})
\left[\begin{array}{c}
0 \\
\bar{n}\bar{c}^T \\
\bar{k}\bar{c}^T
\end{array}\right].
\]
Moreover,
\[\hat{f}_1(\bar{x})=f(\bar{x}(A^{-1})^T)\bar{m}^T=
\bar{x}(A^{-1})^T\bar{m}^T+\alpha=x_1+\alpha.\]
Now our claim follows by the same method as in the proof of Lemma~\ref{staleg}.

\vspace{2ex}

{\em Case } 1c: Suppose that rank $G(\bar{c})=2$.
Then we can assume that $\bar{a},\bar{b},\bar{c}\in{\Bbb{Z}}^3$  and
$\bar{a},\bar{ b}$ generates $G(\bar{c})$. Set $q=\det
A\in{\Bbb{N}}$. Then the linear part of $\hat{f}$ (equals
$K=A\,N\,A^{-1}$) belongs to $M_3({\Bbb{Z}}/q)$.
Moreover,
the functions $\tilde{F}:{\Bbb{R}}^2\rightarrow{\Bbb{R}}^2$,
$\tilde{\gamma}:{\Bbb{R}}^2\rightarrow{\Bbb{R}}$ and
$\tilde{\xi}:{\Bbb{R}}^2\rightarrow{\Bbb{R}}$ are ${\Bbb{Z}}^2$--periodic,
by Corollary~\ref{gener} (see Appendix~\ref{app}).

\vspace{2ex}

{\em Case } 1c(i): Suppose that $d_1/d_2$ is irrational. From
(\ref{syf0}) we obtain $K'=\left[\begin{array}{cc}1&0\\0&1\end{array}\right]$.
Set
\[L:=\left[
\begin{array}{ccc}
1/q & 0 & 0\\
0 & 1/q & 0\\
0 & 0 & 1
\end{array}\right].\]
Now consider the diffeomorphism
$\check{f}:\tor^3_{LA}\rightarrow\tor^3_{LA}$ given by
$\check{f}=L\circ\hat{f}\circ L^{-1}$.
Then
\[\check{f}(x_1,x_2,x_3)=(\check{F}(x_1,x_2),x_3+\check{\gamma}(x_1,x_2)),\]
where
\begin{eqnarray*}
\check{F}(x_1,x_2) & = & \frac{1}{q}F(qx_1,qx_2),\\
\check{\gamma}(x_1,x_2) & = & \gamma(qx_1,qx_2).
\end{eqnarray*}
Then
\begin{eqnarray*}
\check{F}(\bar{x}+\bar{m})-\check{F}(\bar{x}) & = &  \bar{m},\\
\check{\gamma}(\bar{x}+\bar{m})-
\check{\gamma}(\bar{x}) & = &
qK_{31}m_1+qK_{32}m_2\in{\Bbb{Z}}
\end{eqnarray*}
for all $\bar{m}\in {\Bbb{Z}}^2$.
Therefore, $\check{f}$ can be also treated as a diffeomorphism of the torus
$\tor^3$. Let $\check{\xi}(x_1,x_2)=\xi(qx_1,qx_2)$. Then
\begin{equation}\label{kiki}
\check{\xi}\circ\check{F}=\check{\xi}+\alpha,
\end{equation}
$D\check{\xi}:{\Bbb{R}}^2\rightarrow{\Bbb{R}}^2$ is
${\Bbb{Z}}^2$--periodic and non--zero at each point.
Moreover, $\check{f}:\tor^3\rightarrow\tor^3$ has $\tau$--polynomial uniform
growth of the derivative. Precisely,
\begin{equation}\label{polcheck}
\frac{1}{n^{\tau}}D\check{f}^n\rightarrow
\left[\begin{array}{cc}
0 & 0\\
D\check{\xi} & 0
\end{array}
\right]
\end{equation}
uniformly.

Let us denote by $\varphi^t$ the Hamiltonian $C^2$--flow on $\tor^2$ defined by
the Hamiltonian equation
\[\frac{d}{dt}\varphi^t(\bar{x})=
\left[
\begin{array}{r}
\check{\xi}_{x_2}(\varphi^t(\bar{x}))\\-
\check{\xi}_{x_1}(\varphi^t(\bar{x}))
\end{array}
\right].\]
Since $\varphi^t$ has no fixed point and
\begin{eqnarray*}\int_{\tor^2}
\check{\xi}_{x_1}(\bar{x})d\bar{x}/
\int_{\tor^2}
\check{\xi}_{x_2}(\bar{x})d\bar{x}=d_1/d_2
\end{eqnarray*}
is irrational, it follows that $\varphi^t$ is $C^2$--conjugate to the
special flow constructed over the rotation by an irrational number $a$
and under a positive $C^2$--function $b:\tor\rightarrow{\Bbb{R}}$ with
$\int_{\tor}b(x)dx=1$, (see for instance \cite[Ch.\ 16]{Co-Fo-Si}) i.e.\
there exists an area--preserving $C^2$--diffeomorphism
$\rho:{\Bbb{R}}^2\rightarrow{\Bbb{R}}^2$ and a matrix $N\in
GL_2({\Bbb{Z}})$ such that
\[\sigma^t\circ\rho=\rho\circ\varphi^t,\]
where $\sigma^t(x_1,x_2)=(x_1,x_2+t)$ and
\[
\rho(\bar{x}+\bar{m})=(\rho_1(\bar{x})+(\bar{m}N)_1+(\bar{m}N)_2a,\rho_2(\bar{x})-
b^{((\bar{m}N)_2)}(\rho_1(\bar{x})))
\]
for all $\bar{m}\in{\Bbb{Z}}^2$.
Let
$T_{a,-b}:\tor\times{\Bbb{R}}\rightarrow\tor\times{\Bbb{R}}$ denote the skew
product given by $T_{a,-b}(x_1,x_2)=(x_1+a,x_2-b(x_1))$. Let us consider
the quotient space $M=\tor\times{\Bbb{R}}/\sim$, where the relation $\sim$ is
defined by $(x_1,x_2)\sim(y_1,y_2)$  iff  $(x_1,x_2)=T_{a,-b}^k(y_1,y_2)$ for
an integer $k$. Then the quotient flow of the action $\sigma^t$ by the
relation $\sim$ is the special flow constructed over the rotation by $a$ and
under the function $b$. Moreover, $\rho:\tor^2\rightarrow M$ conjugates
flows $\varphi^t$ and $\sigma^t$. Let $\bar{F}:M\rightarrow
M$ denote the $C^2$--diffeomorphism given by
$\bar{F}:=\rho\circ\check{F}\circ\rho^{-1}$.
Since the map ${\Bbb{R}}\ni
t\longmapsto\check{\xi}(\varphi^t\bar{x})\in{\Bbb{R}}$ is constant for
each  $\bar{x}\in{\Bbb{R}}^2$ we see that  the map
\[{\Bbb{R}}\ni t\longmapsto\check{\xi}\circ\rho^{-1}(\sigma^t(x_1,x_2))=
\check{\xi}\circ\rho^{-1}(x_1,x_2+t)\in{\Bbb{R}}\]
is constant for each  $(x_1,x_2)\in{\Bbb{R}}^2$. It follows that the function
$\check{\xi}\circ\rho^{-1}:{\Bbb{R}}^2\rightarrow{\Bbb{R}}$ depends only
on the first coordinate.
Moreover,
\begin{eqnarray*}
D\rho^{-1}(\bar{x})
\left[
\begin{array}{c}
0\\ 1\end{array}
\right] & = & \frac{d}{dt}\rho^{-1}\circ\sigma^t(\bar{x})|_{t=0}=
\frac{d}{dt}\varphi^t\circ\rho^{-1}(\bar{x})|_{t=0}\\
 & = &
\left[
\begin{array}{r}
\check{\xi}_{x_2}(\rho^{-1}(\bar{x}))\\
-
\check{\xi}_{x_1}(\rho^{-1}(\bar{x}))
\end{array}
\right].
\end{eqnarray*}
Consequently, $\partial\rho^{-1}_1/\partial
x_2=\partial\check{\xi}/\partial x_2\circ\rho^{-1}$ and
$\partial\rho^{-1}_2/\partial
x_2=-\partial\check{\xi}/\partial x_1\circ\rho^{-1}$. It follows that
\[\frac{d}{dx_1}\check{\xi}\circ\rho^{-1}=
\frac{\partial \check{\xi}}{\partial x_1}\circ\rho^{-1}\cdot
\frac{\partial \rho^{-1}_1}{\partial x_1}+
\frac{\partial \check{\xi}}{\partial x_2}\circ\rho^{-1}\cdot
\frac{\partial \rho^{-1}_2}{\partial x_1}=-\det D\rho^{-1}=\delta=\pm 1.
\]
Therefore
\begin{equation}\label{dro}
\check{\xi}\circ\rho^{-1}(x_1,x_2)=\delta x_1+c.
\end{equation}
>From (\ref{kiki}) we see that $\check{\xi}\circ\rho^{-1}\circ\bar{F}=
\check{\xi}\circ\rho^{-1}+\alpha$ and consequently that
$\bar{F}_1(x_1,x_2)=x_1+\delta\alpha$. To shorten notation, we will
write $\alpha$ instead of $\delta\alpha$. Since
$\bar{F}:{\Bbb{R}}^2\rightarrow{\Bbb{R}}^2$ preserves area, we conclude that
\[\bar{F}(x_1,x_2)=(x_1+\alpha,\ep x_2+\beta(x_1)),\]
where $\beta:{\Bbb{R}}\rightarrow{\Bbb{R}}$ is a  $C^2$--function and
$\ep=\det D\bar{F}=\pm 1$. As
$\bar{F}$ is a diffeomorphism of $M$, there exists
$m_1,m_2\in{\Bbb{Z}}$ such that
\begin{eqnarray*}
(x_1+1+\alpha,\ep x_2+\beta(x_1+1)) & = & \bar{F}(x_1+1,x_2)=
T_{a,-b}^{m_2}\bar{F}(x_1,x_2)+(m_1,0)\\
& = & (x_1+\alpha+m_1+m_2a,\ep x_2+\beta(x_1)-b^{(m_2)}(x_1+\alpha)).
\end{eqnarray*}
It follows that $m_1=1$, $m_2=0$, hence that
$\beta:\tor\rightarrow{\Bbb{R}}$.
Moreover,  there exists
$n_1,n_2\in{\Bbb{Z}}$ such that
\begin{eqnarray*}
\lefteqn{
(x_1+a+\alpha,\ep x_2-\ep b(x_1)+\beta(x_1+a))}\\
 & = & \bar{F}\circ T_{a,-b}(x_1,x_2)=
T_{a,-b}^{n_2}\bar{F}(x_1,x_2)+(n_1,0)\\
& = & (x_1+\alpha+n_1+n_2a,\ep x_2+\beta(x_1)-b^{(n_2)}(x_1+\alpha)).
\end{eqnarray*}
It follows that $n_1=0$, $n_2=1$, hence that $\beta(x)-b(x+\alpha)=
-\ep b(x)+\beta(x+a)$. Hence
\[1-\ep=\int_{\tor}(b(x+\alpha)-\ep
b(x))dx=\int_{\tor}(\beta(x)-\beta(x+a))dx=0.\]
Therefore
\[\bar{F}(x_1,x_2)=(x_1+\alpha,x_2+\beta(x_1))\]
and the skew products $\bar{F}$ and $T_{a,-b}$ commute.

Let us consider the diffeomorphism
$\bar{f}:M\times\tor\rightarrow M\times\tor$ defined by
\[\bar{f}:=\rho\times\mbox{Id}_{\tor}\circ\check{f}\circ(\rho\times\mbox{
Id}_{\tor})^{-1}.\]
Then
\[\bar{f}(x_1,x_2,x_3)=(\bar{F}(x_1,x_2),x_3+\bar{\gamma}(x_1,x_2)),\]
where $\bar{\gamma}:M\rightarrow\tor$ is given by
$\bar{\gamma}=\check{\gamma}\circ\rho^{-1}$. Therefore there exist
$k_1,k_2\in{\Bbb{Z}}$ such that
\begin{eqnarray*}
\bar{\gamma}(x_1+1,x_2) & = & \bar{\gamma}(x_1,x_2)+k_1\\
\bar{\gamma}(x_1+a,x_2-b(x_1)) & = & \bar{\gamma}(x_1,x_2)+k_2.
\end{eqnarray*}
Moreover,
\begin{eqnarray*}
\frac{1}{n^{\tau}}D\bar{f}^n & = &
\left[\begin{array}{cc}
D\rho\circ\check{F}^n\circ\rho^{-1} & 0\\
0 & 1
\end{array}
\right]
\left[\begin{array}{cc}
n^{-\tau}D\check{F}^n\circ\rho^{-1} & 0\\
n^{-\tau}D(\check{\gamma}^{(n)})\circ\rho^{-1} & n^{-\tau}
\end{array}
\right]
\left[\begin{array}{cc}
D\rho^{-1} & 0\\
0 & 1
\end{array}
\right]
\\
 & \rightarrow &
\left[\begin{array}{cc}
0 & 0\\
D\check{\xi}\circ\rho^{-1} & 0
\end{array}
\right]
\left[\begin{array}{cc}
D\rho^{-1} & 0\\
0 & 1
\end{array}
\right]\\
 &  & =
\left[\begin{array}{cc}
0 & 0\\
D(\check{\xi}\circ\rho^{-1}) & 0
\end{array}
\right]=
\left[\begin{array}{ccc}
0 & 0 & 0\\
0 & 0 & 0\\
\delta & 0 & 0
\end{array}
\right]
\end{eqnarray*}
uniformly, by (\ref{polcheck}) and (\ref{dro}). It follows that
\[\frac{1}{n^{\tau}}\sum_{k=0}^{n-1}(\bar{\gamma}_{x_1}\circ
\bar{F}^k+\bar{\gamma}_{x_2}\circ
\bar{F}^k\cdot D\beta^{(k)})\rightarrow \delta\]
and
\[\frac{1}{n^{\tau}}\sum_{k=0}^{n-1}\bar{\gamma}_{x_2}\circ
\bar{F}^k\rightarrow 0\]
uniformly, hence that
\[\frac{1}{n^{\tau}}\sum_{k=0}^{n-1}\int_{M}(\bar{\gamma}_{x_1}(
\bar{F}^k(x_1,x_2))+\bar{\gamma}_{x_2}(
\bar{F}^k(x_1,x_2))\;D\beta^{(k)}(x_1))dx_1dx_2\rightarrow \delta\]
and
\begin{equation}\label{bec}
\frac{1}{n^{\tau}}\sum_{k=0}^{n-1}\int_{M}\bar{\gamma}_{x_2}
(\bar{F}^k(x_1,x_2))dx_1dx_2\rightarrow 0.
\end{equation}
We show that
\[\frac{1}{n}\sum_{k=0}^{n-1}\int_{M}(\bar{\gamma}_{x_1}(
\bar{F}^k(x_1,x_2))+\bar{\gamma}_{x_2}(
\bar{F}^k(x_1,x_2))D\beta^{(k)}(x_1))dx_1dx_2\rightarrow b(0)k_1.\]
This gives $\tau=1$ and $k_1\neq 0$.
To prove it, note that
\begin{eqnarray*}
\lefteqn{\frac{1}{n}\sum_{k=0}^{n-1}\int_{M}\bar{\gamma}_{x_1}(
\bar{F}^k(x_1,x_2))dx_1dx_2}\\ & = &
\int_0^1\int_0^{b(x_1)}\bar{\gamma}_{x_1}(x_1,x_2)dx_2dx_1\\
& = &
\int_0^1\frac{d}{dx_1}\int_0^{b(x_1)}\bar{\gamma}(x_1,x_2)dx_2dx_1-
\int_0^1Db(x_1)\bar{\gamma}(x_1,b(x_1))dx_1\\
& = &
\int_0^{b(1)}\bar{\gamma}(1,x_2)dx_2-
\int_0^{b(0)}\bar{\gamma}(0,x_2)dx_2-
\int_0^1Db(x_1)(\bar{\gamma}(x_1+a,0)-k_2)dx_1\\
& = &
b(0)k_1-\int_0^1Db(x_1)\bar{\gamma}(x_1+a,0)dx_1.
\end{eqnarray*}
Next observe that
\begin{eqnarray*}
\lefteqn{\frac{1}{n}\sum_{k=0}^{n-1}\int_{M}\bar{\gamma}_{x_2}(
\bar{F}^k(x_1,x_2))D\beta^{(k)}(x_1)dx_1dx_2}\\
& = &
\frac{1}{n}\sum_{k=0}^{n-1}\int_0^1\int_0^{b(x_1)}\bar{\gamma}_{x_2}(x_1,x_2)
D\beta^{(k)}(x_1-k\alpha)dx_2dx_1\\
& = &
\frac{1}{n}\sum_{k=0}^{n-1}\int_0^1(\bar{\gamma}(x_1,b(x_1))-
\bar{\gamma}(x_1,0))
D\beta^{(k)}(x_1-k\alpha)dx_1\\
& = &
\frac{1}{n}\sum_{k=0}^{n-1}\int_0^1(\bar{\gamma}(x_1+a,0)-k_2-
\bar{\gamma}(x_1,0))
D\beta^{(k)}(x_1-k\alpha)dx_1\\
& = &
\frac{1}{n}\sum_{k=0}^{n-1}\int_0^1\bar{\gamma}(x_1+a,0)
(D\beta^{(k)}(x_1-k\alpha)-D\beta^{(k)}(x_1-k\alpha+a))dx_1\\
& = &
\frac{1}{n}\sum_{k=0}^{n-1}\int_0^1\bar{\gamma}(x_1+a,0)
(Db(x_1)-Db(x_1-k\alpha))dx_1\\
& = &
\int_0^1\bar{\gamma}(x_1+a,0)Db(x_1)dx_1-
\int_0^1\bar{\gamma}(x_1+a,0)
\frac{1}{n}\sum_{k=0}^{n-1}Db(x_1-k\alpha)dx_1\\
&  & \rightarrow
\int_0^1\bar{\gamma}(x_1+a,0)Db(x_1)dx_1,
\end{eqnarray*}
which proves the desire conclusion. From (\ref{bec}) we have
\[\int_0^1\int_0^{b(x_1)}\bar{\gamma}_{x_2}
(x_1,x_2)dx_2dx_1=0.\]
However
\begin{eqnarray*}
\int_0^1\int_0^{b(x_1)}\bar{\gamma}_{x_2}
(x_1,x_2)dx_2dx_1
& = &
\int_0^1(\bar{\gamma}(x_1,b(x_1))-\bar{\gamma}(x_1,0))dx_1\\
 & = &
\int_0^1(\bar{\gamma}(x_1+a,0)-\bar{\gamma}(x_1,0)-k_2)dx_1\\
 & = &
k_1a-k_2.
\end{eqnarray*}
It follows that $k_1a=k_2$, which contradicts the fact that $k_1\neq 0$
and $a$ is irrational. Consequently, $d_1/d_2$ must be rational.

\vspace{2ex}

{\em Case } 1c(ii): Suppose that $(d_1,d_2)=d(l_1,l_2)$, where $l_1,l_2$
are relatively
prime integer numbers. Then there exist $M\in GL_2({\Bbb{Z}})$ and
$m\in{\Bbb{Z}}$ such that
\[K'=M^{-1}\left[\begin{array}{cc}
1 & 0 \\ m/q & \ep
\end{array}
\right]M,\]
by (\ref{syf0}). Hence there exists an even number $r$ such that
$K'^r\in GL_2({\Bbb{Z}})$. Therefore  the diffeomorphism
$F^r:{\Bbb{R}}^2\rightarrow{\Bbb{R}}^2$ can  be treated as an
area--preserving diffeomorphism of the torus $\tor^2$.
Let $\check{\xi}:\tor^2\rightarrow\tor$ be given by
$\check{\xi}(x_1,x_2)=d^{-1}\xi(x_1,x_2)$. Then
\begin{eqnarray*}
\check{\xi}\circ F^r=\check{\xi}+r\alpha/d &
\mbox{ and } &
(d_1(\check{\xi}),d_2(\check{\xi}))=(l_1,l_2)\neq 0.
\end{eqnarray*}
Notice that $\alpha/d$ is irrational. Indeed, suppose that
$\alpha/d=k/l$, where $k\in{\Bbb{Z}}$ and $l\in{\Bbb{N}}$. Let
$\Xi:\tor^3_A\rightarrow{\Bbb{C}}$ be given by
\[\Xi(x_1,x_2,x_3)=\exp 2\pi i l\check{\xi}(x_1,x_2).\]
As $\check{\xi}\circ F=\check{\xi}+k/l$, we have
\[\Xi\hat{f}(x_1,x_2,x_3)=\exp 2\pi i
l\check{\xi}\circ F(x_1,x_2)=\Xi(x_1,x_2,x_3).\]
By the ergodicity of $\hat{f}$, $\Xi$ and finally $\check{\xi}$ is
constant, which is impossible.

By Theorem~\ref{conju} (see Appendix~\ref{app0}),
there exists an area--preserving $C^2$--diffeomorphism
$\psi:\tor^2\rightarrow\tor^2$ such that
\[\psi^{-1}\circ F^r\circ\psi:\tor^2\rightarrow\tor^2\]
is a skew product and $\check{\xi}\circ\psi(x_1,x_2)=x_1$. Therefore
$D(\xi\circ\psi)=[\,d\;0\,]$.
Let $L\in GL_2({\Bbb{Z}})$ denote the linear part of $\psi$.
Set
\[\bar{L}:=\left[\begin{array}{cc}L & 0\\ 0& 1\end{array}\right]
\in GL_3({\Bbb{Z}}).\]
Let us consider the area--preserving $C^2$--isomorphism
$\rho:\tor^3_A\rightarrow\tor^3_{\bar{L}^{-1}A}$ defined by
\[\rho(x_1,x_2,x_3)=(\psi^{-1}(x_1,x_2),x_3).\]
Let $\check{f}:\tor^3_{\bar{L}^{-1}A}\rightarrow\tor^3_{\bar{L}^{-1}A}$ be
given by $\check{f}=\rho\circ\hat{f}\circ\rho^{-1}$.
Then
\begin{eqnarray*}
\frac{1}{n^{\tau}}D\check{f}^n & = &
\left[\begin{array}{cc}
D\psi^{-1}\circ F^n\circ\psi & 0\\
0 & 1
\end{array}
\right]
\left[\begin{array}{cc}
n^{-\tau}DF^n\circ\psi & 0\\
n^{-\tau}D(\gamma^{(n)})\circ\psi & n^{-\tau}
\end{array}
\right]
\left[\begin{array}{cc}
D\psi & 0\\
0 & 1
\end{array}
\right]
\\
 & \rightarrow &
\left[\begin{array}{cc}
0 & 0\\
D\xi\circ\psi & 0
\end{array}
\right]
\left[\begin{array}{cc}
D\psi & 0\\
0 & 1
\end{array}
\right]\\
 &  & =
\left[\begin{array}{cc}
0 & 0\\
D(\xi\circ\psi) & 0
\end{array}
\right]=
\left[\begin{array}{ccc}
0 & 0 & 0\\
0 & 0 & 0\\
d & 0 & 0
\end{array}
\right]
\end{eqnarray*}
uniformly. Let us consider the diffeomorphism
$\bar{f}:\tor^3\rightarrow\tor^3$ given by
$\bar{f}=A^{-1}\bar{L}\circ\check{f}\circ\bar{L}^{-1}A$.
It is easy to see that
\[\frac{1}{n^{\tau}}D\bar{f}^n\rightarrow
A^{-1}\bar{L}
\left[\begin{array}{ccc}
0 & 0 & 0\\
0 & 0 & 0\\
d & 0 & 0
\end{array}
\right]
\bar{L}^{-1}A\]
uniformly and that $\bar{f}$ and $f$ are conjugate via the
area--preserving
$C^2$--diffeomorphism $A^{-1}\bar{L}\circ\rho\circ A:\tor^3\rightarrow\tor^3$.
Now applying Lemma~\ref{staleg} for $\bar{f}$ gives our claim.

\vspace{2ex}

{\em Case} 2: Suppose that $g=(h_1\bar{a}^T+h_2\bar{b}^T)\bar{c}$, where
$\bar{a}$ and $\bar{b}$ are orthogonal to $\bar{c}$ and the determinant of
the matrix
\[
A^{-1}=\left[\begin{array}{ccc}
\bar{c}^T &
\bar{a}^T &
\bar{b}^T
\end{array}\right]\]
equals $1$.
Let $\hat{f}:\tor^3_A\rightarrow\tor^3_A$ is given by $\hat{f}:=A\circ
f\circ A^{-1}$.
Then
\[\hat{g}=
A(h_1\bar{a}^T+h_2\bar{b}^T)\bar{c}
A^{-1}
=\left[\begin{array}{c}
0\\
\hat{h}_1 \\
\hat{h}_2
\end{array}\right]
\left[\begin{array}{ccc}
1 &
0 &
0
\end{array}\right],\]
where $\hat{h}_i(\bar{x}):=h_i(A^{-1}\bar{x})$ for $i=1,2$.
>From (\ref{kapel}) we obtain
\[\left[\begin{array}{c}
0\\
\hat{h}_1(\bar{x}) \\
\hat{h}_2(\bar{x})
\end{array}\right]
\left[\begin{array}{ccc}
1 &
0 &
0
\end{array}\right]  =
\left[\begin{array}{c}
0\\
\hat{h}_1(\hat{f}\bar{x}) \\
\hat{h}_2(\hat{f}\bar{x})
\end{array}\right]
\left[\begin{array}{ccc}
1 &
0 &
0
\end{array}\right]D\hat{f}(\bar{x})
\]
for all $\bar{x}\in\tor^3_A$.
Consequently
\[\frac{\partial}{\partial x_1}\hat{f}_1(\bar{x})\hat{h}_i(\hat{f}\bar{x})=
\hat{h}_i(\bar{x}),\;\;\;
\frac{\partial}{\partial x_2}\hat{f}_1(\bar{x})=
\frac{\partial}{\partial x_3}\hat{f}_1(\bar{x})=0\;\;
\mbox{ and }\;\;
\frac{\partial}{\partial x_1}\hat{f}_1(\bar{x})\neq 0\]
for all $\bar{x}\in\tor^3_A$ and $i=1,2$.
Now observe that $\hat{h}_1,\hat{h}_2$ are linearly
dependent. Indeed, without loss of generality we can assume that
$\hat{h}_2$ is $A\lambda^{\otimes3}$ non--zero. Then
$\hat{h}_2(\bar{x})\neq 0$ for a.e.\
$\bar{x}\in\tor^3_A$, by the ergodicity of $\hat{f}$.
Therefore the measurable function
$\hat{h}_1/\hat{h}_2:\tor^3_A\rightarrow{\Bbb{R}}$ is
$\hat{f}$--invariant. Hence there is a real constant $c$ such that
$\hat{h}_1(\bar{x})=c\hat{h}_2(\bar{x})$ for a.e.\
$\bar{x}\in\tor^3_A$, by ergodicity. Consequently, ${h}_1=c{h}_2$, which
leads us to Case 1, and the proof is complete. $\Box$

\section{$4$--dimensional case}\indent\label{4d}

In this section we indicate why there is no $4$--dimensional analogue of
classifications of area--preserving diffeomorphisms of polynomial growth
of the derivative presented in previous sections.
Precisely, we construct an ergodic area--preserving diffeomorphism of the
$4$--dimensional torus with linear uniform growth of the derivative which
is not even metrically isomorphic to any $3$--steps skew product, i.e.\ to any
automorphism of $\tor^4$ of the form
\[\tor^4\ni (x_1,x_2,x_3,x_4)\longmapsto
(x_1+\alpha,\ep_1x_2+\beta(x_1),\ep_2x_3+\gamma(x_1,x_2),\ep_3x_4+\delta(x_1,x_2,x_3))\in\tor^4,\]
where $\ep_i=\pm 1$ for $i=1,2,3$.
Before we pass to the $4$--dimensional case we should mention area
preserving diffeomorphisms of the $2$--torus with sublinear growth of the
derivative.

\begin{df}
{\em We say that a $C^1$--diffeomorphism $f:\tor^2\rightarrow\tor^2$ has}
sublinear growth of the derivative {\em if the sequence $Df^n/n$ tends
uniformly to zero.}
\end{df}

Suppose that $f:\tor^2\rightarrow\tor^2$ is an area--preserving weakly mixing
$C^{\infty}$--diffeomorphism with sublinear growth of the derivative. The
examples of such diffeomorphisms will be given later. Let
$T_{\varphi}:\tor^2\rightarrow\tor^2$ be an Anzai skew product where
$Tx=x+\alpha$ is an ergodic rotation on the circle and
$\varphi:\tor\rightarrow\tor$ is a $C^{\infty}$--function with non--zero
topological degree.

\begin{theo}
The product diffeomorphism
$f\times T_{\varphi}:\tor^4\rightarrow\tor^4$ is ergodic and has linear
uniform growth of the derivative. Moreover, it is not metrically isomorphic
to any $3$--steps skew product.
\end{theo}

\pf The first hypothesis of the theorem is obvious.
Now suppose, contrary to our claim,
that $f\times T_{\varphi}$ is metrically isomorphic to a $3$--steps skew
product. Then $f\times T_{\varphi}$ is measure theoretically distal (has
generalized discrete spectrum in the terminology of \cite{Zi}).
However, $f\times T_{\varphi}$ has a weakly mixing factor,
which is impossible, because measure theoretically distal and weakly
mixing dynamical systems are disjoint (see~\cite{Fu}). $\Box$

\vspace{2ex}

In the remainder of this section we present two examples of
area--preserving weakly mixing
diffeomorphisms with sublinear growth of the derivative.

For given $\alpha\in\tor$ and $\beta:\tor\rightarrow{\Bbb{R}}$ we will
denote by
$T_{\alpha,\beta}:\tor\times{\Bbb{R}}\rightarrow\tor\times{\Bbb{R}}$ the
skew product $T_{\alpha,\beta}(x_1,x_2)=(x_1+\alpha,x_2+\beta(x_1))$.
Let $a\in\tor$ be an irrational number and let
$b:\tor\rightarrow{\Bbb{R}}$ be a positive $C^{\infty}$--function with
$\int_{\tor}b(x)dx=1$. Let $\sigma^t$ denote the flow on
$\tor\times{\Bbb{R}}$ given by $\sigma^t(x_1,x_2)=(x_1,x_2+t)$.
Let us consider the quotient space $M_{a,b}=\tor\times{\Bbb{R}}/\sim$,
where the relation $\sim$ is defined by $(x_1,x_2)\sim(y_1,y_2)$  iff
$(x_1,x_2)=T_{a,-b}^k(y_1,y_2)$ for an integer $k$.
Then the quotient flow $\sigma_{a,b}^t$ of the action $\sigma^t$ by the
relation $\sim$ is the special flow constructed over the rotation by $a$ and
under the function $b$.
By Lemma 2 in \cite{Fa} and Theorem 1 in \cite{Mo}, there exists a
$C^{\infty}$--diffeomorphism $\rho:M_{a,b}\rightarrow\tor^2$ such that the
flow $\varphi^t=\rho\circ\sigma_{a,b}^t\circ\rho^{-1}$ is a Hamiltonian
flow on $\tor^2$ with no fixed points, i.e.\ there exists
$C^{\infty}$--function $\xi:{\Bbb{R}}^2\rightarrow{\Bbb{R}}$
such that $D\xi$ is ${\Bbb{Z}}^2$--periodic, non--zero at each point and
\[\frac{d}{dt}\varphi^t(\bar{x})=
\left[
\begin{array}{r}
{\xi}_{x_2}(\varphi^t(\bar{x}))\\-
{\xi}_{x_1}(\varphi^t(\bar{x}))
\end{array}
\right].\]
We will identify $\rho$ with a diffeomorphism
$\rho:{\Bbb{R}}^2\rightarrow{\Bbb{R}}^2$ such that
\begin{eqnarray*}
\rho(x_1+1,x_2) & = & \rho(x_1,x_2)+(N_{11},N_{12}),\\
\rho(x_1+a,x_2-b(x_1)) & = & \rho(x_1,x_2)+(N_{21},N_{22})
\end{eqnarray*}
for any $(x_1,x_2)\in{\Bbb{R}}^2$, where $N\in GL_2({\Bbb{Z}})$.
Then
\begin{eqnarray}
\label{rho1}
D\rho(x_1+1,x_2) & = & D\rho(x_1,x_2) \\
\label{rho2}
D\rho(T_{a,-b}^n(x_1,x_2))
\left[\begin{array}{cc}
1 & 0 \\ -Db^{(n)}(x_1) & 1
\end{array}\right] & = & D\rho(x_1,x_2)
\end{eqnarray}
for any integer $n$.

Let $T_{\alpha,\beta}:\tor\times{\Bbb{R}}\rightarrow\tor\times{\Bbb{R}}$
be a commuting with $T_{a,-b}$ skew product, where
$\beta:\tor\rightarrow{\Bbb{R}}$ is of class $C^{\infty}$. Then
$T_{\alpha,\beta}$ can be treated as the $C^{\infty}$--diffeomorphism of
$M_{a,b}$. Let us denote by $f:\tor^2\rightarrow\tor^2$ the area--preserving
$C^{\infty}$--diffeomorphism  given by $f:=\rho\circ
T_{\alpha,\beta}\circ \rho^{-1}$.

\begin{lem}\label{sub}
The diffeomorphism  $f:\tor^2\rightarrow\tor^2$ has sublinear growth of
the derivative.
\end{lem}

\pf Since
\[Df^n(\bar{x})=D\rho(T_{\alpha,\beta}^n\circ \rho^{-1}(\bar{x}))
\left[\begin{array}{cc}
1 & 0 \\ D\beta^{(n)}(\rho_1^{-1}(\bar{x})) & 1
\end{array}\right] D\rho^{-1}(\bar{x}),\]
it suffices to show that
\[\frac{1}{n}D\rho(T_{\alpha,\beta}^n(x_1,x_2))
\left[\begin{array}{cc}
1 & 0 \\ D\beta^{(n)}(x_1) & 1
\end{array}\right] \rightarrow 0\]
uniformly on the set $M'=\{(x_1,x_2):x_1\in{\Bbb{R}},0\leq x_2\leq b(x_1)\}$.
For every $(x_1,x_2)\in{\Bbb{R}}^2$ let us denote by $n(x_1,x_2)$ the unique
integer number such that $T_{a,-b}^{n(x_1,x_2)}(x_1,x_2)\in M'$, i.e.\
\[b^{(n(x_1,x_2))}(x_1)\leq x_2\leq b^{(n(x_1,x_2)+1)}(x_1).\]
Let $c,C$ be positive constants such that $0<c\leq b(x)\leq C$ for every
$x\in\tor$. Then
\[ c|n(x_1,x_2)|\leq |x_2|\leq C|n(x_1,x_2)|+C.\]
Since
\begin{eqnarray*}
\lefteqn{\frac{1}{n}D\rho(T_{\alpha,\beta}^n(x_1,x_2))
\left[\begin{array}{cc}
1 & 0 \\ D\beta^{(n)}(x_1) & 1
\end{array}\right]}\\
 & = &
D\rho(T_{a,-b}^{n(T_{\alpha,\beta}^n(x_1,x_2))}T_{\alpha,\beta}^n(x_1,x_2))
\frac{1}{n}\left[\begin{array}{cc}
1 & 0 \\ -Db^{(n(T_{\alpha,\beta}^n(x_1,x_2)))}(x_1+n\alpha)+
D\beta^{(n)}(x_1) & 1
\end{array}\right]
\end{eqnarray*}
(by~(\ref{rho2})), $D\rho$ is bounded on $M'$ (by~(\ref{rho1})) and
$n^{-1}D\beta^{(n)}$ tends uniformly to zero, it suffices to
show that
\[\frac{1}{n}Db^{(n(T_{\alpha,\beta}^n(x_1,x_2)))}(x_1+n\alpha)\rightarrow
0\] uniformly on $M'$. To prove it, first observe that
\[|n(T_{\alpha,\beta}^n(x_1,x_2))|\leq c^{-1}|x_2+\beta^{(n)}(x_1)|\leq
k_1+k_2n,\]
for any natural $n$ and every $(x_1,x_2)\in M'$, where $k_1=C/c$ and
$k_2=\|\beta\|_{\sup}/c$. Fix $\ep>0$.
Let $n_0$ be a natural number such that $|n|\geq n_0$ implies
\begin{eqnarray*}
\frac{1}{|n|}\|Db^{(n)}\|_{\sup}<\ep/2k_2 & \mbox{ and } & k_1+k_2n\leq 2k_2
\end{eqnarray*}
for any integer $n$. Assume that $n$ is a natural number such that $n\geq
\|b\|_{C^1}n_0/\ep$. Let $(x_1,x_2)\in M'$.
If $|n(T_{\alpha,\beta}^n(x_1,x_2))|\|b\|_{C^1}/n<\ep$, then
\[|\frac{1}{n}Db^{(n(T_{\alpha,\beta}^n(x_1,x_2)))}(x_1+n\alpha)|\leq
\frac{|n(T_{\alpha,\beta}^n(x_1,x_2))|}{n}\|b\|_{C^1}<\ep.\]
Otherwise,
\[|n(T_{\alpha,\beta}^n(x_1,x_2))|\geq \ep n/\|b\|_{C^1}\geq n_0.\]
Then
\begin{eqnarray*}
\lefteqn{|\frac{1}{n}Db^{(n(T_{\alpha,\beta}^n(x_1,x_2)))}(x_1+n\alpha)|}\\
& \leq &
\frac{|n(T_{\alpha,\beta}^n(x_1,x_2))|}{n}
\frac{1}{|n(T_{\alpha,\beta}^n(x_1,x_2))|}
\|Db^{(n(T_{\alpha,\beta}^n(x_1,x_2)))}\|_{\sup}\\
 & < &
\frac{k_1+k_2n}{n}\ep/2k_2\leq \ep,
\end{eqnarray*}
which completes the proof. $\Box$

\begin{prop} (see \cite{Aar})\label{p1}
For every $C^2$--function $\beta:\tor\rightarrow{\Bbb{R}}$ with
zero mean, which is not a trigonometric polynomial there exists a dense
$G_{\delta}$ set of irrational numbers $\alpha\in\tor$ such that the
corresponding skew product
$T_{\alpha,\beta}:\tor\times{\Bbb{R}}\rightarrow\tor\times{\Bbb{R}}$ is
ergodic.
\end{prop}

>From the proof of the Main Theorem in \cite{Sz} and the nature of the weak
mixing property, one can obtain the following

\begin{prop} \label{p2}
For every positive real analytic function $b:\tor\rightarrow{\Bbb{R}}$ with
$\int_{\tor}b(x)dx=1$, which is not a trigonometric polynomial there exists
a dense $G_{\delta}$ set of irrational numbers $a\in\tor$ such that the
corresponding special flow $\sigma_{a,b}^t$ is weakly mixing.
\end{prop}

{\bf Example 1.} Suppose that $\sigma_{a,b}^t$ is a weakly mixing special
flow whose roof function is real analytic. Let $\varphi^t$ be a Hamiltonian
flow on $\tor^2$ that is $C^{\infty}$--conjugate to the special flow
$\sigma_{a,b}^t$. Then the area--preserving diffeomorphism
$\varphi^1:\tor^2\rightarrow\tor^2$ is weakly mixing and has sublinear
growth of the derivative, by Lemma~\ref{sub}.

\vspace{2ex}

{\bf Example 2.} By Propositions~\ref{p1} and \ref{p2},
there exist a $C^{\infty}$--function $\beta:\tor\rightarrow{\Bbb{R}}$ with
zero mean and an irrational numbers $\alpha\in\tor$ such that the
corresponding skew product
$T_{\alpha,\beta}:\tor\times{\Bbb{R}}\rightarrow\tor\times{\Bbb{R}}$ is
ergodic and for no real $r\neq 0$ there exist $c\in\tor$ and a measurable
function $c_r:\tor\rightarrow\tor$ satisfying
\[c_r(x+\alpha)e^{2\pi ir\beta(x)}=c\cdot c_r(x).\]
Assume that $b:\tor\rightarrow{\Bbb{R}}$ is a positive $C^{\infty}$--function
with $\int_{\tor}b(x)dx=1$ and $a\in\tor$ is an irrational numbers such that
the skew products $T_{\alpha,\beta}$  and $T_{a,-b}$ commute and
$T_{a,-b+1}$ is not measurably isomorphic to any iteration of
$T_{\alpha,\beta}$. Let us consider the special flow $\sigma_{a,b}^t$ on
$M_{a,b}$ and the diffeomorphism $T_{\alpha,\beta}:M_{a,b}\rightarrow
M_{a,b}$. Then  $T_{\alpha,\beta}:M_{a,b}\rightarrow
M_{a,b}$ is ergodic, by the ergodicity of
$T_{\alpha,\beta}:\tor\times{\Bbb{R}}\rightarrow\tor\times{\Bbb{R}}$.
Moreover, $T_{\alpha,\beta}:M_{a,b}\rightarrow M_{a,b}$ is weakly mixing.
Indeed, suppose, contrary to our claim, that  $T_{\alpha,\beta}$ has an
eigenvalue $c\in\tor$. Then there exists a measurable function
$F:\tor\times{\Bbb{R}}\rightarrow\tor$ such that
\begin{eqnarray}\label{ccc}
F\circ T_{\alpha,\beta}=cF & \mbox{ and } & F\circ T_{a,-b}=F.
\end{eqnarray}
It follows that
\[F\circ\sigma^t\circ T_{\alpha,\beta}=F\circ T_{\alpha,\beta}\circ\sigma^t=
cF\circ\sigma^t,\]
and hence that
\[F\sigma^t/F\circ T_{\alpha,\beta}=F\sigma^t/F\]
for any real $t$. By the ergodicity of $T_{\alpha,\beta}$, there exists a
measurable function $c:{\Bbb{R}}\rightarrow\tor$ such that $F\sigma^t=c(t)F$.
Then
\[c(t+s)F=F\sigma^{t+s}=F\sigma^{t}\sigma^{s}=c(t)c(s)F\]
for all real $t,s$. Therefore there exists a real number $r$ such
that $c(t)=e^{2\pi irt}$. Consequently,
\[F(x_1,x_2)=F(x_1)e^{2\pi irx_2},\]
where $F:\tor\rightarrow\tor$ is given by $F(x)=F(x,0)$.
>From (\ref{ccc}), we conclude that
\[F(x_1+\alpha)e^{2\pi ir(x_2+\beta(x_1))}=c
F(x_1)e^{2\pi irx_2},\]
and finally that
\[F(x_1+\alpha)e^{2\pi ir\beta(x_1)}=c
F(x_1),\]
which contradicts our assumption.

Let $\varphi^t$ be a Hamiltonian flow on $\tor^2$ that is
$C^{\infty}$--conjugate to the special flow $\sigma_{a,b}^t$, via a
$C^{\infty}$--diffeomorphism $\rho:M_{a,b}\rightarrow\tor^2$. Let us
consider the area--preserving diffeomorphism $f:\tor^2\rightarrow\tor^2$
given by  $f:=\rho\circ T_{\alpha,\beta}\circ \rho^{-1}$. Then $f$ is
weakly mixing and has sublinear
growth of the derivative, by Lemma~\ref{sub}.

\appendix

\section{}\indent\label{app0}

For given $\alpha\in{\Bbb{R}}$, $\varphi:\tor\rightarrow\tor$ and $\ep=\pm 1$ let
$T_{\alpha,\varphi,\ep}:\tor^2\rightarrow\tor^2$ denote the diffeomorphism
$T_{\alpha,\varphi,\ep}(x_1,x_2)=(x_1+\alpha,\ep x_2+\varphi(x_1))$.

\begin{theo}\label{conju}
Let $f:\tor^2\rightarrow\tor^2$ be an
area--preserving $C^2$--diffeomorphism. Suppose that there exist an irrational
number $\alpha$ and a $C^2$--function
$\xi:\tor^2\rightarrow\tor$
such that
\begin{itemize}
\item $\xi\circ f=\xi +\alpha$,
\item $0<c\leq\|D\xi(\bar{x})\|$ for any $\bar{x}\in\tor^2$,
\item $(d_1(\xi),d_2(\xi))=(p_1,p_2)\neq (0,0)$, where $p_1,p_2$ are
relatively prime integer numbers.
\end{itemize}
Then there exist an area--preserving $C^2$--diffeomorphism
$\psi:\tor^2\rightarrow\tor^2$ and a $C^2$--cocycle
$\varphi:\tor\rightarrow\tor$ such that
\[ f\circ \psi =\psi\circ T_{\alpha,\varphi,\ep}\]
where $\ep=\det Df$ and $\xi\circ\psi(x_1,x_2)=x_1$.
\end{theo}

\pf
Without restriction of generality we can assume that $\xi(0,0)=0$. For
every $s\in{\Bbb R}$ set
\[A_s=\{\bar{x}\in{\Bbb R}^2:\;\xi(\bar{x})=s\}.\]
Since $0<c\leq\|D\xi(\bar{x})\|$ for any $\bar{x}\in{\Bbb R}^2$, we
see that $A_s$ is a curve. Moreover,
\begin{equation}\label{ogran}
(x_1,x_2)\in A_s\Rightarrow (x_1+p_2,x_2-p_1)\in A_s\;\;\mbox{ and }\;\;
\bar{x}\in A_s\Rightarrow f(\bar{x})\in A_{s+\alpha}.
\end{equation}
Let $q_1,q_2$ be integer numbers such that $p_1q_1+p_2q_2=1$.
Then
\[(x_1,x_2)\in A_s\Rightarrow (x_1+q_1,x_2+q_2)\in A_{s+1}.\]
Let $\gamma:{\Bbb R}\rightarrow {\Bbb R}^2$ be a
$C^2$--function such that:
\begin{itemize}
\item $\gamma(0)=(0,0)$,
\item $\gamma(s+1)=\gamma(s)+(q_1,q_2)$
for any real $s$,
\item $\xi\circ \gamma(s)=s$ for any real $s$  (see Fig.\ 1).
\end{itemize}
\begin{figure}[here]
\begin{center}
\unitlength 1cm
\begin{picture}(10,8)(-2,-2)
\put(2,-1.75){ Figure 1.}

\put(-1.5,0){\vector(1,0){8}}
\put(0,-1){\vector(0,1){7.5}}
\put(0,0){\circle*{0.075}}
\put(-0.8,0.2){$\scriptstyle\gamma(0)=(0,0)$}
\put(2,5){\circle*{0.075}}
\put(1.35,5.15){$\scriptstyle(p_{\scriptscriptstyle 2},-p_{\scriptscriptstyle 1})$}

\bezier{500}(-1,-1)(-1.75,0)(0,0)
\bezier{500}(0,0)(1,0)(1.25,1)
\bezier{500}(1.25,1)(1.75,3)(1,4)
\bezier{500}(1,4)(0.25,5)(2,5)
\bezier{500}(2,5)(3,5)(3.25,6)
\put(3,1){\circle*{0.075}}
\put(1.9,0.6){$\scriptstyle\gamma(1)=
(q_{\scriptscriptstyle 1},q_{\scriptscriptstyle 2})
$}
\put(5,6){\circle*{0.075}}
\put(4.2,5.7){$\scriptstyle(q_{\scriptscriptstyle 1}+p_{\scriptscriptstyle 2},
q_{\scriptscriptstyle 1}-p_{\scriptscriptstyle 1})$}
\put(0.2,2.5){$A_0$}
\put(0.7,2.65){\vector(1,0){0.7}}

\bezier{500}(2.5,-1)(2.5,-0.5)(2,0)
\bezier{500}(2,0)(1.25,1)(3,1)
\bezier{500}(3,1)(4,1)(4.25,2)
\bezier{500}(4.25,2)(4.75,4)(4,5)
\bezier{500}(4,5)(3.25,6)(5,6)
\put(4.7,5){$A_1$}
\put(4.7,5.15){\vector(-1,0){0.7}}

\bezier{500}(-1,0)(-0.5,0.25)(0,0)
\bezier{500}(0,0)(1,-0.5)(1.5,0.5)
\bezier{500}(1.5,0.5)(2,1.5)(3,1)
\bezier{500}(3,1)(4,0.5)(4.5,1.5)
\put(3.4,4){$A_s$}
\put(3.4,4.15){\vector(-1,0){0.5}}

\bezier{500}(0.5,-1)(1.5,-0.5)(1.5,1)
\bezier{500}(1.5,1)(1.5,1.5)(2,2)
\bezier{500}(2,2)(2.5,2.5)(2,3)
\bezier{500}(2,3)(1.5,3.5)(2.5,4)
\bezier{500}(2.5,4)(3.5,4.5)(3.5,6)

\put(2.4,2){trace of $\gamma$}
\put(2.7,1.9){\vector(-1,-2){0.35}}

\put(1.47,0.44){\circle*{0.075}}
\put(1.6,-0.7){$\gamma(s)$}
\put(1.75,-0.35){\vector(-1,3){0.24}}

\end{picture}
\end{center}
\end{figure}

For every $s\in{\Bbb R}$, denote by
$\psi(s,\,\cdot\,):{\Bbb R}\rightarrow{\Bbb R}^2$ the solution of the
following differential equation
\[\begin{array}{rcl}
\frac{d}{dt}\psi(s,t) & = &
\left[\begin{array}{r}
-\xi_{x_2}(\psi(s,t))\\
\xi_{x_1}(\psi(s,t))
\end{array}\right]\\
\psi(s,0) & = & \gamma(s).
\end{array}\]
Clearly, $\psi:{\Bbb R}^2\rightarrow{\Bbb R}^2$ is of class $C^2$ and
\[\psi(s+1,t)=\psi(s,t)+(q_1,q_2)\]
for any $(s,t)\in{\Bbb R}^2$. Moreover,
\[\frac{d}{dt}\xi\circ \psi(s,t)=
D\xi(\psi(s,t))\cdot \frac{d}{dt}\psi(s,t)=0.\]
Hence
\begin{equation}\label{xis}
\xi(\psi(s,t))=\xi(\psi(s,0))=\xi(\gamma(s))=s.
\end{equation}
It follows that
\begin{equation}\label{curv}
A_s=\{\psi(s,t):t\in{\Bbb R}\},
\end{equation}
because $\|\frac{d}{dt}\psi(s,t)\|\geq c$. Therefore $\psi$ is a surjection.
Moreover,
\begin{eqnarray}\label{prut}
\det D\psi(s,t) & = &
\det
\left[\begin{array}{cr}
\frac{d}{ds}\psi_1(s,t) & -\xi_{x_2}(\psi(s,t))\\
\frac{d}{ds}\psi_2(s,t) & \xi_{x_1}(\psi(s,t))
\end{array}\right]\\
 & = & D\xi(\psi(s,t))\frac{d}{ds}\psi(s,t)=\frac{d}{ds}\xi\circ\psi(s,t)
=1.\nonumber
\end{eqnarray}
Next note that $\psi$ is an injection. Suppose that $(s_1,t_1)\neq
(s_2,t_2)$ and $\psi(s_1,t_1)=\psi(s_2,t_2)$. Then
\[s_1=\xi\circ\psi(s_1,t_1)=\xi\circ\psi(s_2,t_2)=s_2,\]
by (\ref{xis}). Hence $t_1\neq t_2$ and $\psi(s_1,t_1)=\psi(s_1,t_2)$.
It follows that  $\psi(s_1,t+t_1-t_2)=\psi(s_1,t)$ for any real $t$, by
the definition of $\psi$. Therefore, the curve $A_{s_1}$ is close,
contrary to (\ref{ogran}).
Hence $\psi$ is an area--preserving $C^2$--diffeomorphism of ${\Bbb
R}^2$.

Since
\[\xi(f\circ\psi(s,t))=\xi(\psi(s,t))+\alpha=s+\alpha,\]
(\ref{curv}) shows that $f\circ\psi(s,t)=\psi(s+\alpha,\eta(s,t))$, where
$\eta:{\Bbb R}^2\rightarrow{\Bbb R}$ is a $C^2$--function.
By assumption, $D\xi\circ f\cdot Df=D\xi$.
It follows that
\[\left[\begin{array}{r}
-\xi_{x_2}\circ f\\
\xi_{x_1}\circ f
\end{array}\right]=\ep Df\left[\begin{array}{r}
-\xi_{x_2}\\
\xi_{x_1}
\end{array}\right],\]
Therefore
\begin{eqnarray*}
\frac{d}{dt}f\circ\psi(s,t)  & = &
Df(\psi(s,t))
\frac{d}{dt}\psi(s,t)
 =
Df(\psi(s,t))
\left[\begin{array}{r}
-\xi_{x_2}(\psi(s,t))\\
\xi_{x_1}(\psi(s,t))
\end{array}\right]\\
& = &
\ep\left[\begin{array}{r}
-\xi_{x_2}(f\circ\psi(s,t))\\
\xi_{x_1}(f\circ\psi(s,t))
\end{array}\right]
 =
\ep\left[\begin{array}{r}
-\xi_{x_2}(\psi(s+\alpha,\eta(s,t)))\\
\xi_{x_1}(\psi(s+\alpha,\eta(s,t)))
\end{array}\right]\\
& =  &
\ep\psi_{x_2}(s+\alpha,\eta(s,t)).
\end{eqnarray*}
On the other hand
\[\frac{d}{dt}f\circ\psi(s,t)=
\frac{d}{dt}\psi(s+\alpha,\eta(s,t))=
\psi_{x_2}(s+\alpha,\eta(s,t))\frac{d}{dt}\eta(s,t).\]
Hence $\eta(s,t)=\ep t+\eta(s,0)$. Let
$\varphi:{\Bbb R}\rightarrow{\Bbb R}$ be given by
$\varphi(s)=\eta(s,0)$. Then
\[f\circ\psi(s,t)=\psi(s+\alpha,\ep t+\varphi(s)).\]

We only need to show that $\psi$ is a diffeomorphism of the torus $\tor^2$.
>From (\ref{ogran}) and (\ref{curv}), for every real $s$ there exists a
unique real number $\tau(s)\neq 0$ such that
\[\psi(s,\tau(s))=\psi(s,0)+(p_2,-p_1).\]
Then the function $\tau:{\Bbb R}\rightarrow{\Bbb R}$ is continuous,
$\tau(s+1)=\tau(s)$ and
\[\psi(s,t+\tau(s))=\psi(s,t)+(p_2,-p_1)\]
for all real $s,t$, by the definition of $\psi$.
Let $M\in GL_2({\Bbb Z})$ denote the linear part of $f$.
As $\xi\circ f=\xi+\alpha$, we have
\[\left[\,p_1\;p_2\right]M=
\left[\,p_1\;p_2\right].\]
It follows that
\[M
\left[\begin{array}{r}
p_2\\-p_1
\end{array}\right]
=\overline{\ep}
\left[\begin{array}{r}
p_2\\-p_1
\end{array}\right],\]
where $\overline{\ep}=\det M=\pm 1$.
Hence
\[f((x_1,x_2)+(p_2,-p_1))=f(x_1,x_2)+\overline{\ep}(p_2,-p_1)\]
for any $(x_1,x_2)\in{\Bbb R}^2$. Therefore
\begin{eqnarray*}
\psi(s+\alpha,\tau(s)+\varphi(s)) & = & f\circ\psi(s,\ep\tau(s))
=  f(\psi(s,0)+\ep (p_2,-p_1))\\
& = &  f\circ\psi(s,0)+\ep\overline{\ep}(p_2,-p_1)
=\psi(s+\alpha,\varphi(s))+\ep\overline{\ep}(p_2,-p_1)\\
 & = & \psi(s+\alpha,\ep\overline{\ep}\tau(s+\alpha)+\varphi(s)).
\end{eqnarray*}
Therefore $\tau(s+2\alpha)=\tau(s)$ for any real $s$. Since $\alpha$ is
irrational, $\tau$ is continuous and periodic of period $1$,
we see that $\tau$ is constant. Hence
\[\psi(s+1,t)=\psi(s,t)+(q_1,q_2),\]
\[\psi(s,t+\tau)=\psi(s,t)+(p_2,-p_1)\]
for any $(s,t)\in{\Bbb R}^2$.
Let $\psi_\tau:{\Bbb R}^2\rightarrow{\Bbb R}^2$ be given by
$\psi_\tau(s,t)=\psi(s,\tau t)$. Then $\psi_{\tau}$ is a
diffeomorphism of $\tor^2$. We conclude from (\ref{prut}) that $\det
D\psi_\tau(s,t)=\tau$ for
all $(s,t)\in{\Bbb R}^2$, hence that $\tau=\pm 1$. It follows that
\[\psi(s+1,t)=\psi(s,t)+(q_1,q_2),\]
\[\psi(s,t+1)=\psi(s,t)\pm (p_2,-p_1)\]
and finally that
$\psi$ is a $C^2$--diffeomorphism of $\tor^2$ such that
\[f\circ\psi=\psi\circ T_{\alpha,\varphi,\ep}.\;\;\;\Box\]

\section{}\indent\label{app}

Let $g:{\Bbb{R}}^2\rightarrow{\Bbb{R}}$ be a continuous function and let
$\bar{a},\bar{b}\in{\Bbb{R}}^3$ be linearly independent vectors. Let
$\bar{c}\in{\Bbb{R}}^3$ be a non--zero vector orthogonal to both $\bar{a}$
and $\bar{b}$.

\begin{lem}\label{per}
Suppose that there exists a vector $\bar{d}\in{\Bbb{R}}^3$ such that
\[h(x_1+\bar{a}\bar{m}^T,x_2+\bar{b}\bar{m}^T)=h(x_1,x_2)+\bar{d}\bar{m}^T\]
for all $\bar{m}\in{\Bbb{Z}}^3$. Then there exist $k_1,k_2\in{\Bbb{R}}$ such
that $\bar{d}=k_1\bar{a}+k_2\bar{b}$ and the function
$\tilde{h}(x_1,x_2)=h(x_1,x_2)-k_1x_1-k_2x_2$ is
$(\bar{a}\bar{m}^T,\bar{b}\bar{m}^T)$--periodic for all
$\bar{m}\in{\Bbb{Z}}^3$. Moreover,
\begin{itemize}
\item if {\em rank} $G(\bar{c})$=0, then $\tilde{h}$ is constant;
\item if {\em rank} $G(\bar{c})$=1, then there exit $l_1,l_2\in{\Bbb{R}}$
and a continuous function $\rho:\tor\rightarrow {\Bbb{R}}$ such that
$\tilde{h}(x_1,x_2)=\rho(l_1x_1+l_2x_2)$ and
$l_1\bar{a}+l_2\bar{b}\in{\Bbb{Z}}^3$ generates $G(\bar{c})$;
\item if {\em rank} $G(\bar{c})$=2, then $\bar{c}\in c{\Bbb{Z}}^3$ where
$c\neq 0$.
\end{itemize}
\end{lem}

\pf Let us consider the function $H:{\Bbb{R}}^3\rightarrow{\Bbb{R}}$ given
by $H(\bar{x}):=h(\bar{a}\bar{x}^T,\bar{b}\bar{x}^T)$ for all
$\bar{x}\in{\Bbb{R}}^3$. Then
$H(\bar{x}+\bar{m})=H(\bar{x})+\bar{d}\bar{m}^T$ for all
$\bar{m}\in{\Bbb{Z}}^3$ and
\[H((x_1,x_2,x_3)(A^{-1})^T)=h(\bar{a}A^{-1}\bar{x}^T,\bar{b}A^{-1}\bar{x}^T)=h(x_1,x_2),\]
where
\[
A=\left[\begin{array}{c}
\bar{a} \\
\bar{b} \\
\bar{c}
\end{array}\right].\]
It follows that the function
$\tilde{H}:{\Bbb{R}}^3\rightarrow{\Bbb{R}}$ given by
$\tilde{H}(\bar{x})=H(\bar{x})-\bar{d}\bar{x}^T$ is
${\Bbb{Z}}^3$--periodic. Let
$A^{-1}=\left[\,\bar{a}'^T\;\bar{b}'^T\;\bar{c}'^T\,\right]$.
Then
\[\bar{d}\bar{a}'^Tx_1+\bar{d}\bar{b}'^Tx_2+\bar{d}\bar{c}'^Tx_3=
h(x_1,x_2)-\tilde{H}((x_1,x_2,x_3)(A^{-1})^T).\]
As $\tilde{H}$ is bounded, we have $\bar{d}\bar{c}'^T=0$.
Define
\[\tilde{h}(x_1,x_2):=h(x_1,x_2)-k_1x_1-k_2x_2,\]
where $k_1=\bar{d}\bar{a}'^T$ and $k_2=\bar{d}\bar{b}'^T$. Then
$k_1\bar{a}+k_2\bar{b}=\bar{d}$ and
\[
\tilde{h}(x_1,x_2)=\tilde{H}((x_1,x_2,x_3)(A^{-1})^T)\]
for all $(x_1,x_2,x_3)\in {\Bbb{R}}^3$. Therefore
\begin{eqnarray*}
\tilde{h}(x_1+\bar{a}\bar{m}^T,x_2+\bar{b}\bar{m}^T) & = &
\tilde{H}(((x_1,x_2,x_3)+\bar{m}A^T)(A^{-1})^T)\\
& = &
\tilde{H}(((x_1,x_2,x_3)(A^{-1})^T+\bar{m})=
\tilde{h}(x_1,x_2)
\end{eqnarray*}
for all $\bar{m}\in{\Bbb{Z}}^3$. Moreover,
\[\tilde{H}(x_1,x_2,x_3)=\tilde{h}(\bar{a}\bar{x}^T,\bar{b}\bar{x}^T).\]
The function $\tilde{H}$ can be  represented as the Fourier series
\[\tilde{H}(\bar{x})=\sum_{\bar{m}\in{\Bbb{Z}}^3}a_{\bar{m}}\exp 2\pi
i\bar{m}\bar{x}^T,\]
where the series converges in $L^2(\tor^3)$.
It follows that
\[\tilde{h}(x_1,x_2)=\sum_{\bar{m}\in{\Bbb{Z}}^3}a_{\bar{m}}\exp 2\pi
i(\bar{m}\bar{a}'^Tx_1+\bar{m}\bar{b}'^Tx_2+\bar{m}\bar{c}'^Tx_3)\]
in $L^2(\tor^3_A)$. Therefore, $\bar{m}\bar{c}'^T=0$, wherever
$a_{\bar{m}}\neq 0$. Hence
\[\tilde{H}(\bar{x})=\sum_{\bar{m}\in G(\bar{c})}a_{\bar{m}}\exp 2\pi
i\bar{m}\bar{x}^T,\]
because $G(\bar{c})=G(\bar{c}')$.

Suppose that rank $G(\bar{c})=0$. Then $\tilde{H}$, and consequently
$\tilde{h}$ is constant.

Now suppose that rank $G(\bar{c})=1$. Then
$\tilde{H}(\bar{x})=\rho(\bar{m}\bar{x}^T)$, where $\bar{m}$ is a
generator of $G(\bar{c})$ and $\rho:\tor\rightarrow{\Bbb{R}}$ is a
continuous function. Moreover,
\[\tilde{h}(x_1,x_2)=\rho(\bar{m}A^{-1}\bar{x}^T)=
\rho(\bar{m}\bar{a}'^Tx_1+\bar{m}\bar{b}'^Tx_2)\]
and
\[\bar{m}\bar{a}'^T\bar{a}+\bar{m}\bar{b}'^T\bar{b}=\bar{m}A^{-1}A=\bar{m}.\]

Now assume that rank $G(\bar{c})=2$. Suppose that
$\bar{m},\bar{n}\in{\Bbb{Z}}^3$ generate $G(\bar{c})$. Then $\bar{c}=c\,\bar{m}\times\bar{n}\in
c\,{\Bbb{Z}}^3$. $\Box$

\begin{lem}
Let $\bar{c}\in{\Bbb{Z}}^3\setminus\{0\}$.
Then there exists a pair of generators $\bar{a},\bar{b}\in{\Bbb{Z}}^3$ of
$G(\bar{c})$ such that
$\Lambda(\bar{a},\bar{b})=\{(\bar{a}\bar{m}^T,\bar{b}\bar{m}^T)\in{\Bbb{Z}}^2
:\bar{m}\in{\Bbb{Z}}^3\}={\Bbb{Z}}^2$.
\end{lem}

\pf Our proof starts with the observation that we can assume that
$c_1,c_2,c_3$ are relatively prime. Let
$p_1:=\gcd(c_2,c_3)$, $p_2:=\gcd(c_3,c_1)$, $p_3:=\gcd(c_1,c_2)$. Then
$p_1,p_2,p_3$ are pairwise relatively prime and
there exist pairwise relatively prime integer numbers  $k_1,k_2,k_3$ such
that $c_1=k_1p_2p_3$, $c_2=p_1k_2p_3$, $c_1=p_1p_2k_3$.
Let us choose integer numbers $x,y$ such that $k_2x+k_3y=1$. Then
\begin{eqnarray*}
\bar{a}:=(p_1,-p_2k_1x,-p_3k_1y) &\mbox{ and } & \bar{b}:=(0,-p_2k_3,p_3k_2)
\end{eqnarray*}
generates the group $G(\bar{c})$. Indeed, suppose that $(m_1,m_2,m_3)\in
G(\bar{c})$. Since $m_1c_1=-m_2c_2-m_3c_3$ we see that $p_1|m_1c_1$,
hence that $p_1|m_1$, because $c_1$ and $p_1$ are relatively prime.
Similarly, there exist $n_1,n_2,n_3\in{\Bbb{Z}}$ such that $m_1=p_1n_1$,
 $m_2=p_2n_2$, $m_3=p_3n_3$. As $(m_1,m_2,m_3)\in G(\bar{c})$
we have
\[n_1k_1+n_2k_2+n_3k_3=0.\]
Since $n_1k_1+n_1(-k_1x)k_2+n_1(-k_1y)k_3=0$ we obtain
\[(n_2+n_1k_1x)k_2+(n_3+n_1k_1y)k_3=0.\]
Therefore there exists $s\in{\Bbb{Z}}$ such that $n_2+n_1k_1x=-sk_3$ and
$n_3+n_1k_1y=sk_2$. It follows that
\[(n_1,n_2,n_3)=n_1(1,-k_1x,-k_1y)+s(0,-k_3,k_2),\]
hence that $(m_1,m_2,m_3)=n_1\bar{a}+s\bar{b}$.
Let $r,t$ be integer numbers such that $rp_1-tc_1=1$. Then
\[(\bar{a},\bar{b})(r,tp_3k_2,tp_2k_3)^T=(1,0).\]
Moreover
\[(\bar{a},\bar{b})(0,r',t')^T=(s',1),\]
whenever $-p_2k_3r'+p_3k_2t'=1$. It follows that $(1,0)$ and $(s',1)$
belong to $\Lambda(\bar{a},\bar{b})$, and finally that
$\Lambda(\bar{a},\bar{b})={\Bbb{Z}}^2$. $\Box$

\begin{cor}\label{gener}
For any pair of generators $\bar{a}',\bar{b}'\in{\Bbb{Z}}^3$ of
$G(\bar{c})$ we have
$\Lambda(\bar{a}',\bar{b}')={\Bbb{Z}}^2$.
\end{cor}

\pf
Suppose that $\bar{a}',\bar{b}'\in{\Bbb{Z}}^3$  generate
$G(\bar{c})$.
Then there exists $K\in GL_2({\Bbb{Z}})$ such that
$(\bar{a}'^T,\bar{b}'^T)=(\bar{a}^T,\bar{b}^T)K$. Therefore
\[\Lambda(\bar{a}',\bar{b}')=\Lambda(\bar{a},\bar{b})K={\Bbb{Z}}^2.\;\;\;
\Box\]

\noindent
Faculty of Mathematics and Computer Science,\\
Nicholas Copernicus University\\
ul. Chopina 12/18\\
87-100 Toru\'n, Poland

\vspace{2ex}\noindent
Institute of Mathematics\\
Polish Academy of Science\\
ul. Chopina 12\\
87-100 Toru\'n, Poland

\vspace{2ex}\noindent
E-mail: fraczek@mat.uni.torun.pl

 \end{document}